%% file: main.tex
\documentclass[12pt]{amsart}
\usepackage{amsmath}
\usepackage{amssymb}
\usepackage{amsthm}
\usepackage{epsfig,graphicx}
\usepackage{mathrsfs}

\newcommand{\ignore}[1]{\relax}

\newcommand{\Id}{\operatorname{Id}}

\newcommand{\Def}{\operatorname{Def}}

\newcommand{\C}{\mathbb C}

\newcommand{\R}{\mathbb R}
\newcommand{\Z}{\mathbb Z}
\newcommand{\N}{\mathbb N}

\newcommand{\DD}{\mathcal D}

\newtheorem{lem}{Lemma}[section]
\newtheorem{claim}{Claim}[section]

\newtheorem{theorem}{Theorem}
\newtheorem{corollary}[lem]{Corollary}
\newtheorem{thm}[lem]{Theorem}

\newtheorem{coro}[lem]{Corollary}
\newtheorem{prop}[lem]{Proposition}

\theoremstyle{definition}
\newtheorem{condition}[claim]{Condition}
\newtheorem{defn}[lem]{Definition}

\newtheorem{problem}[claim]{Problem}
\newtheorem{exa}[lem]{Example}
\newtheorem{que}[lem]{Question}

\theoremstyle{remark}
\newtheorem{rmk}[lem]{Remark}

\newtheorem{ack}{Acknowledgment}

\renewcommand{\setminus}{\smallsetminus}

\newcommand{\ctor}{(\C^\times)^{n}}
\newcommand{\tordva}{(\C^\times)^{2}}

\newcommand{\ev}{\operatorname{ev}}

\newcommand{\dd}{\partial}

\newcommand{\cp}{{\mathbb C}{\mathbb P}}
\newcommand{\rp}{{\mathbb R}{\mathbb P}}

\newcommand{\Log}{\operatorname{Log}}

\newcommand{\Int}{\operatorname{Int}}

\renewcommand{\setminus}{\smallsetminus}

\newcommand{\MM}{{\mathcal M}}

\newcommand{\vol}{\operatorname{vol}}

\newcommand{\Hom}{\operatorname{Hom}}

\ignore{
\newtheorem{theorem}{Theorem}

\newtheorem{condition}[theorem]{Condition}

\newenvironment{proof}[1][Proof]{\noindent\textbf{#1.} }{\ \rule{0.5em}{0.5em}}
}

\renewcommand{\-}{\hspace{-3pt}}
\renewcommand{\setminus}{\smallsetminus}
\newcommand{\mdn}{\mu_\Delta\circ\nu_\epsilon}
\newcommand{\tme}{\tilde\mu_\epsilon}
\newcommand{\sm}{\smallsetminus}

\begin{document}
\title
{Examples of tropical-to-Lagrangian correspondence}
\author{Grigory Mikhalkin}
\address{Universit\'e de Gen\`eve,  Math\'ematiques, Battelle Villa, 1227 Carouge, Suisse}
\begin{abstract}
The paper associates 
Lagrangian submanifolds in symplectic
toric varieties
to certain tropical curves
inside the convex polyhedral domains of $\R^n$
that appear as the images of the moment map of the toric varieties.

We pay a particular attention to the case $n=2$,
where we reprove Givental's theorem \cite{Gi}
on Lagrangian embeddability
of non-oriented surfaces to $\C^2$,
as well as to the case $n=3$, where we
see appearance of the graph 3-manifolds
studied by Waldhausen \cite{Wald} as Lagrangian
submanifolds.
In particular, rational tropical curves in $\R^3$
produce 
3-dimensional rational homology spheres.
The order of their first homology groups
is determined by the
multiplicity of tropical curves in the corresponding
enumerative problems.
\end{abstract}
\thanks{Research is supported in part by the grants 159240, 159581 and the NCCR SwissMAP
project of the Swiss National Science Foundation.}
\maketitle

\section{Some background material}
\subsection{Symplectic toric varieties and the moment
map.}
\input stv.tex
\subsection{Tropical curves in $\Delta$}
\input trg.tex

\input results.tex

\input exa2.tex

\section{Proof of Theorem \ref{main}} 
\label{s-proof}
\input proof.tex

\input more-bkg.tex

\section{Rational tropical curves and three-dimensional Lagrangians}
\input l3.tex

\begin{ack}
The author had benefited from
many useful discussions with
Tobias Ekholm, Yakov Eliashberg,
Sergey Galkin, Alexander Givental,
Ilia Itenberg,
Conan Leung, Diego Matessi, Vivek Shende
and Oleg Viro.
\end{ack}

\bibliography{b}
\bibliographystyle{plain}

\end{document}

%% file: stv.tex
Let $\Lambda\approx\Z^n$ be a free Abelian
group (a lattice) of rank $n$. 
Let $A\approx\R^n$ be an affine space 
over the real $n$-dimensional
vector space $\Lambda\otimes\R$.
Clearly we can identify
$T_xA=\Lambda\otimes\R$ for the tangent
space to $A$ at any $x\in A$. In particular,
the tangent spaces at all points of $A$ are
canonically identified.

Furthermore, as $A$ is an affine space over $T_xA$
a choice of $x$ gives an identification
between $A$ and $\Lambda\otimes\R$,
$A\ni y\mapsto y-x\in \Lambda\otimes\R$.
Thus an element $p$ of the dual lattice
$\Lambda^*=\Hom(\Lambda,\Z)$
and a point $x\in A$ define an affine function
$p^x:A\to \R$ by $y\mapsto\ <\- p,y-x\- >$.

We refer to $A$ as the {\em tropical affine space}
(more classically it is an affine space corresponding
to the structure group obtained from the integer 
linear group $GL(n,\Z)$ by extending it with all
real translations in $\R^n$).

\begin{defn}
A {\em polyhedral domain} $\Delta\subset A$
is the
intersection of a finite number of half-spaces
$\{x\in A\ |\ p^y_jx\ge a_j\}$, $j=1,\dots,N$, $y\in A$,
$a_j\in \R$ and $p_j\in\Lambda^*$ such that
the interior of $\Delta$ is non-empty.

A point $z\in\Delta$ is called {\em smooth}
if there exists an open neighborhood $U\ni z$,
a point $y\in \R^n$ and an integer basis
$\{p_1,\dots,p_n\}\subset\Lambda^*$
such that
$$U\cap\Delta=\bigcap\limits_{j=1}^n
\{x\in\R^n\ |\ p^y_jx\ge 0\}.$$
A polyhedral domain is called
{\em Delzant},
if all its points are smooth.
\end{defn}

Suppose that $(M,\omega)$ is a $2n$-dimensional
symplectic manifold with a Hamiltonian action
of the real $n$-torus $T=(S^1)^n$, and 
$\mu:M\to t^*$ is the corresponding moment map,
defined by 
\begin{equation}\label{mudef}
<\- Y,d\mu(X)\- >\ =2\pi\omega(X,Y),
\end{equation}
for any $Y\in t$, $X\in T_uM$, $u\in M$.
Here $t^*$ is
the dual vector space to the Lie 
algebra $t$ of $T$.
The element $Y\in t$ yields a vector field on $M$
through the action of $T$ on $M$, thus
the left-hand side of \eqref{mudef} is well-defined for
any $X\in T_uM$.
The condition \eqref{mudef} defines $\mu$ up to
a translation in $t^*$,
see e.g. \cite{CdS} for details.

Note that the Lie algebra $t$ comes with a natural
lattice $\Lambda^*\approx\Z^n$ defined
as $\exp^{-1}(0)$ for the exponent map $\exp:t\to T$
so that $T=t/\Lambda^*$.
We have
$t^*=\Lambda\otimes\R$ for
$\Lambda=\Hom(\Lambda^*,\Z)$.
Since $\mu$ is defined up to a translation,
its target has a natural structure of an affine
space over $A$, i.e. the tropical affine space.
We write
\begin{equation}
\mu:M\to A
\end{equation}
for the moment map of $(M,\omega)$.
The fibers of the moment map coincide with the
orbits of the action of $T$.

\begin{defn}
A {\em symplectic toric variety} is a $2n$-dimensional
symplectic manifold $(M,\omega)$
with a Hamiltonian action
of the real $n$-torus $T$ such that $\mu(M)$ is
a polyhedral domain.
\end{defn}

\begin{rmk}\label{nonconvD}
If $(M,\omega)$ is compact $2n$-dimensional
symplectic manifold with a Hamiltonian action
of the real $n$-torus $T$ then $\mu(M)$ is
automatically a polyhedral domain, and furthermore
is always smooth. For non-compact $M$ then
the condition that $\mu(M)$ is polyhedral is 
a certain completeness condition on $\omega$.
\end{rmk}

For any Delzant polyhedral domain $\Delta\subset A$
there exists a symplectic toric variety
$(M_\Delta,\omega_\Delta)$
with $\mu(M_\Delta)=\Delta$, see \cite{De}.
\begin{exa}
Suppose $\Delta=A$. The cotangent space $T^*A$
has a canonical symplectic structure $\omega=d\alpha$,
$\alpha=pdq$, $q\in A$, $p\in T_q A$ (so that
$\alpha$ is a well-defined non-closed 1-form).
The Lie algebra $t=\Lambda^*\otimes\R$ acts
on $T^*A$ by translations preserving
the form $\omega$. Thus the quotient space
$M=(T^*A)/\Lambda^*$ together with the form $\omega$
is a symplectic manifold
with an action of $T=t/\Lambda^*$.
According to \eqref{mudef}
the moment map is given by the projection onto $A$,
\[
\mu:M\to A,\ (p,q)\mapsto q.
\]

The quotient
$(M,\omega)$
may be identified with the complex $n$-torus
$\ctor$
enhanced with a $\ctor$-invariant symplectic form
\[
\omega_{\ctor}=\frac i2\sum\limits_{j=1}^n 
\frac{dz_j}{z_j}\wedge\frac{d\bar z_j}{\bar z_j}.
\]
The torus $T=(S^1)^n$ acts on $\ctor$ by coordinatewise
multiplication (we identify $S^1$
with a unit circle in $\C$).
The moment map $\mu:\ctor\to A$ 
coincides with $\Log:\ctor\to\R^n$,
\[
\Log(z_1,\dots,z_n)=(\log|z_1|,\dots,\log|z_n|)
\]
after the identification of $A$ with $\R^2$.
In other words, we have a symplectomorphism
between 
$(\ctor,\omega_{\ctor})$ and
$(M,\omega=\sum\limits_{j=1}^n
dp_j\wedge dq_j)$ given by
\[
q_j=\log|z_j|,\
p_j=\arg(z_j).
\]
\end{exa}

If $\Delta\neq A$ then for each (open) face
$E\subset\dd\Delta$ of codimension $k$
the inverse image $\mu^{-1}(E)$ is fibered
by $k$-dimensional tori whose tangent vectors
are contained in the radical of the form
$\omega$ restricted to the tangent space of $\mu^{-1}(E)$.
Taking the quotient of $\mu^{-1}(\Delta)$
by this fibration over all the
faces of $\dd\Delta$ is known as the symplectic 
reduction. In the case when $\Delta$ is a Delzant
polyhedral domain, this construction produces
a smooth $2n$-dimensional symplectic manifold
$(M_\Delta=\mu^{-1}(\Delta)/\sim,\omega_\Delta)$
with a
Hamiltonian action of $T$ such that the
moment map $\mu$
descends by the projection $\mu^{-1}(\Delta)\to M_\Delta$
to the moment map $\mu_\Delta:M_\delta\to A$,
and the symplectic forms $\omega_\Delta$ and $\omega$
agree on
$M_\Delta\setminus\mu^{-1}(\dd\Delta)\subset M$.
In particular, we have $\mu_\Delta(M_\Delta)=\Delta$.

%% file: trg.tex
Let $\bar\Gamma$ be a topological space homeomorphic
to a finite graph. The subset
$\dd\Gamma\subset\bar\Gamma$
of 1-valent vertices and the subset $V_\Gamma$
of vertices of valence greater than 2 do not
depend on the choice of a graph model for $\bar\Gamma$,
and thus are well-defined subspaces in the topological
space $\bar\Gamma$.

We set $\Gamma=\bar\Gamma\setminus\dd\Gamma$.
Connected components
of $\Gamma\setminus V_\Gamma$
are homeomorphic to open intervals and
are called {\em edges} of $\Gamma$.
The set of edges of $\Gamma$ is denoted
with $E_\Gamma$.
An edge is called a {\em leaf} if it is 
adjacent to $\dd\Gamma$ and and a {\em bounded edge}
otherwise.

\begin{defn}
\label{def-trc}
The topological space $\Gamma$ enhanced
with an inner complete metric
is called a (smooth, irreducible
and explicit) {\em tropical curve}.
\end{defn}
Specifying an inner metric on
$\Gamma$ amounts 
to specifying positive lengths of the edges
of the graph. By the completeness assumption
the lengths of the leaves must be infinite
while the lengths of the bounded edges are finite.
\begin{rmk}
Definition \ref{def-trc} introduces
smooth, irreducible and explicit tropical curves.
For the purposes of this paper we refer
to such curves simply as {\em tropical curves}.
In a more general framework, there is a {\em genus
function} (cf. e.g. \cite{Ca})
$V_\Gamma\to\Z_{\ge 0}$
which can also be reformulated as a $\chi$-measure
(cf. e.g. \cite{KaMi}). Then the set $V_\Gamma$ also
includes vertices of valence 1 or 2 with positive
values of the genus function. 
Such generalized curves are needed for
compactifications of moduli spaces of tropical
curves. In this paper we do not use them.
Our tropical curves have the genus function equal
to zero everywhere on $V_\Gamma$.
\end{rmk}
Recall that a continuous map $f:X\to Y$
between topological spaces is called
an {\em immersion} if it is a local embedding,
i.e. for any $x\in X$ there exists an open neighborhood
$U\ni x$
such that $f|_U$ is an embedding of $U$ to $Y$.

\begin{defn}
\label{tr-imm}
An immersion $h:\Gamma\to A$ (between topological
spaces $\Gamma$ and $A$) is called 
{\em tropical} \cite{Mi05},
if the following conditions hold.
\begin{itemize}
\item[(a)] For any edge $e\in E_\Gamma$,
a point $x\in e$ and a unit tangent vector $u\in T_xe$
the restriction $h|_e$ is a smooth map such that
\[
dh_x (u)\in\Lambda\subset T_{h(x)}A.
\]
In particular, $dh_x(u)$ does not depend
on the choice of $x\in e$ and depends only
on the orientation of $e$ given by $u$.
We denote $dh(e)=dh_x(u)\in\Lambda$
for the oriented edge $e$.
\item[(b)]
For any vertex $v\in V_\Gamma$ we have
\[
\sum\limits_{\bar e\ni v} dh(e) = 0,
\]
where the sum is taken over all edges $e$ adjacent
to $v$. The orientation of $e$ is chosen to be
away from $v$. This condition is known as
the {\em balancing condition}.
\end{itemize}
The tropical immersion is called {\em locally flat}
if, in addition, the following condition holds.
\begin{itemize}
\item[(c)] If a collection of edges
$\{e_j\}\subset E_\Gamma$
is adjacent to the same vertex $v\in V_\Gamma$
then the linear span of $\{dh(e_j)\}$ in 
$\Lambda\otimes\R$ is 2-dimensional.
\end{itemize}
The intersection $h(\Gamma)\cap\Delta\subset \Delta$
is called
{\em a locally flat tropical curve in
a polyhedral domain $\Delta$}. 
\end{defn}
Note that if $\Gamma$ is 3-valent then any
tropical immersion $h:\Gamma\to A$ is
locally flat as a consequence of the balancing condition.

\begin{defn}\label{pti}
A tropical immersion
$h:\Gamma\to A$ is called {\em primitive} if
the following
conditions hold.
\begin{itemize}
\item[(i)] The tropical curve $\Gamma$ is connected,
3-valent, and $V_\Gamma\neq\emptyset$.
\item[(ii)] For any $e\in E_\Gamma$ the vector $dh(e)$
is primitive element of the lattice $\Lambda$.
\item[(iii)] If $x\neq y\in\Gamma$ and 
$h(x)=h(y)$ then $x,y\in\Gamma\setminus V_\Gamma$.
\end{itemize}
The image $h(\Gamma)$ is called
a {\em primitive tropical curve in $A$}.
\end{defn}

\begin{defn}
Let $\Delta\subset A$ be a polyhedral domain
and $\Gamma_\Delta$ be a topological space
homeomorphic to a connected graph.
An immersion $h_\Delta:\Gamma_\Delta\to\Delta$
is called {\em $\Delta$-tropical} 
(or just tropical)
if there exists
a tropical immersion $h:\Gamma\to A$ such that
$\Gamma_\Delta\subset\Gamma$,
$h|_{\Gamma_\Delta}=h_\Delta$ and 
$h^{-1}(\dd\Delta)$ is a finite set disjoint from
$V_\Gamma$. A $\Delta$-tropical immersion $h_\Delta$
is called {\em primitive} if $h$ can be chosen
to be primitive and $\#(h^{-1}(x))=1$ whenever
$x\in\dd\Delta$. The symbol $\#$ stands
for the cardinality of a set.

A subset $C\subset\Delta$ is called 
a {\em primitive tropical curve in $\Delta$}
if there exists a primitive $\Delta$-tropical
immersion $h_\Delta:\Gamma_\Delta\to \Delta$
such that $h_\Delta(\Gamma_\Delta)=C$.
\end{defn}

\begin{prop}\label{Sigma-finite}
If $C\subset \Delta$ is a primitive tropical curve
in $\Delta$ then
a primitive $\Delta$-tropical immersion
with $C=h_\Delta(\Gamma_\Delta)\cap\Delta$ is unique.
Furthermore, the
self-intersection set 
\begin{equation}\label{def-Sigma}
\Sigma(C)=\{z\in C\ |\
\#(h_\Delta^{-1}(z))>1\}
\end{equation}
of $h_\Delta$ is finite.
\end{prop}
\begin{proof}
Suppose that $h_\Delta:\Gamma_\Delta\to \Delta$,
$h'_{\Delta}:\Gamma'_\Delta\to \Delta$
are two primitive $\Delta$-tropical immersions and
$h:\Gamma\to A$, $h':\Gamma'\to A$ be tropical
immersions with $h_\Delta=h|_{\Gamma_\Delta}$
and $h'_\Delta=h'|_{\Gamma'_\Delta}$.
Let $v\in V_\Gamma\cap \Gamma_\Delta$.
By the conditions $(i)$ and $(ii)$ of 
Definition \ref{pti},
and since $h$ is a topological immersion,
a small neighborhood of $v$ in $\Gamma_\Delta$
and that of $h(v)$ in $C$
are homeomorphic. Thus $h(v)=h'(v')$ for
$v'\in V_{\Gamma'}\cap\Gamma'_\Delta$. 
By the condition $(iii)$,
we get a 1-1 correspondence between
$V_\Gamma\cap\Gamma_\Delta$ and
$V_{\Gamma'}\cap\Gamma'_\Delta$.
Similarly, we get a 1-1 correspondence
between $h^{-1}(\dd\Delta)$ and $h'^{-1}(\dd\Delta)$.
Also by $(iii)$
if $e_1,e_2\in E_\Gamma$ are such that
$h(e_1)\neq h(e_2)$ then $h(e_1)$ is 
transversal to $h(e_2)$.
This implies that the obtained correspondences
extend to 
a homeomorphism $\Phi:\Gamma_\Delta\to\Gamma'_\Delta$
such that $h'_\Delta=h_\Delta\circ\Phi$.
The same property implies
the finiteness of $\Sigma$.
\end{proof}
Thus we may speak of the vertices
$V_C=h_{\Delta}(\Gamma_\Delta)$ of $C$.

The {\em boundary points} of
a primitive tropical
curve $C$ in $\Delta$
are the points of its (topological) boundary 
$\dd C=C\cap\dd\Delta$.
Each boundary point belongs to a unique
$(n-k)$-dimensional (relatively open) face
of the polyhedral domain $\Delta$. 
We call $k$ the {\em codimension} of a boundary
point $x\in\dd C$.

\begin{defn}
Let $x\in \dd C$ be a boundary point of codimension 1,
$e_x\in E_\Gamma$ be the edge containing $h^{-1}(x)$,
and $\Delta_x\subset\dd\Delta$ be the facet containing
$x$.
The {\em boundary momentum} $p(x)\in\N$
is the {\em tropical intersection number} of
$h(e_x)$ and $\Delta_x$ in $A$, i.e. the index
in $\Lambda$ of 
the sublattice generated by
$dh(e_x)$ and the elements of $\Lambda$
parallel to $\Delta_x$.
\end{defn}

If $y\in\dd C$ is a point of codimension $2$ then
it belongs to the closure of exactly two facets
$\Delta_1$ and $\Delta_2$. As above, we may define
the boundary momentum of $y$ with respect to $\Delta_j$,
$j=1,2$, as the tropical intersection number of 
$h(e_y)$ and $\Delta_j$.
\begin{defn}
A boundary point $y\in\dd C$ of codimension 2
is called a {\em bissectrice} if both of
these boundary momenta are equal to 1.
\end{defn}

\begin{defn}
A primitive tropical curve in $\Delta$ 
is called {\em even} if
\begin{enumerate}
\item
all of its boundary points
are of codimension at most 2, 
\item 
all of its codimension 1 boundary points have
boundary momenta equal to 2 and
\item
all of its codimension 2 boundary points are
bissectrice points.
\end{enumerate}
\end{defn}

Let $v\in V_C$ be a vertex of a primitive 
curve $C$ in $\Delta$ adjacent to 
the edges $h(e_j)$, $j=1,2,3$,
$e_j\in E_\Gamma$ where $h:\Gamma\to A$
is a tropical immersion
with $C=h(\Gamma)\cap\Delta$.
\begin{defn}
The {\em multiplicity} of $v$ is the number
\begin{equation}\label{mult-v}
m(v)=|dh(e_1)\wedge dh(e_2)|,
\end{equation}
i.e. the area of the parallelogram spanned by
the vectors $dh(e_1),dh(e_2)\in\Lambda$.
\end{defn}
Because of the balancing condition 
of Definition \ref{tr-imm} we have
$m(v)=|dh(e_j)\wedge dh(e_k)|$ for any
$j\neq k=1,2,3$.

The {\em self-intersection number} of $v$
is defined as 
\[
\delta(v)=\frac{m(v)-1}2.
\]
Let $w\in\Sigma(C)$ (which is a finite
set by Proposition \ref{Sigma-finite}).
The {\em multiplicity} of $w$ is the number 
\begin{equation}\label{mult-w}
m(w)=\sum\limits_{x\neq y\in h^{-1}(w)}
|dh(e_x)\wedge dh(e_y)|.
\end{equation}
\begin{defn}
The {\em self-intersection number}
of a primitive tropical curve $C$ in $\Delta$
is defined as
\[
\delta(C)=\sum\limits_{v\in V_C}\delta(v)+
\sum\limits_{w\in\Sigma(C)}m(w).
\]
The curve $C$ is called {\em smooth}
if $\delta(C)=0$.
\end{defn}


%% file: results.tex
\section{Main result}
Let $\Delta\subset A$
be a Delzant polyhedral domain,
$M_\Delta$ be the symplectic toric variety
corresponding to $\Delta$, and
$\mu_\Delta:M_\Delta\to\Delta$ be the moment map.
Let $C$
be a primitive tropical curve in $\Delta$,
and $V_C$ be the set of its vertices.

For a vertex $v\in V_C$ we denote with $A_v\subset A$
the 2-dimensional affine subspace containing
the edges adjacent to $v$.
Primitivity of $v$ implies that there are three such
edges, and that these edges do not overlap, so
such $A_v$ is unique.
Similarly, for an edge $e\in E_C$
we denote with $A_e\subset A$
the 1-dimensional affine subspace containing $e$.

Consider a metric on $A\approx\R^n$ invariant
under translations in $A$.
Denote with
$U_\epsilon(v)\subset A_v$ the
intersection of $A_v$ and 
the open ball of radius $\epsilon$ around $v$.
Denote $U_\epsilon(V_C)=
\bigcup\limits_{v\in V_C}U_\epsilon(v)$

Recall that (in 2D-topology)
a {\em pair-of-pants} $P$ is a smooth
surface diffeomorphic to a thrice punctured sphere.
We denote by $P_\delta$ the pair-of-pants with $\delta$
nodes, i.e. the topological space obtained
from $P$ by gluing $\delta$ disjoint pairs of points
to $\delta$ nodes of $P_\Delta$.

Let $A_f\subset A$ be a $k$-dimensional affine subspace
of $A$ with an integer slope, i.e. such that
its tangent space $TA_f\subset TA=\Lambda\otimes \R$
is generated by
a $k$-dimensional sublattice of $\Lambda$.
This is equivalent to requiring that the
{\em conormal space} $N^*A_f\subset T^*A=
\Lambda^*\otimes\R$ is generated by an $(n-k)$-dimensional
sublattice of $\Lambda^*$.

The fiber torus
$\Theta=(\Lambda^*\otimes\R)/\Lambda^*\approx
(S^1)^n$
of the moment map $\mu_\Delta$
can be considered as a (commutative and compact) 
Lie group.
Denote with $\Theta_f\subset\Theta$
the $(n-k)$-dimensional subtorus of $\Theta$ obtained
as $N^*A_f/(\Lambda^*\cap N^*A_f)$.
If $y\in M_{\Int\Delta}=M_\Delta\setminus
\mu^{-1}_\Delta(\dd\Delta)$ then $y+\Theta_f$
is an affine subtorus of the (torus) fiber of 
$\mu_\Delta:M_\Delta\to\Delta$ containing $y$.
Here we are using the sum notations as we have
the action of the Abelian group $\Theta$
on $M_{\Int \Delta}$ coming from 
the action of $\Lambda^*\otimes\R$ on $T^*A$.

Recall that $\nu:L\to M_\Delta$ is called
a {\em Lagrangian immersion} if
$L$ is a smooth $n$-dimensional manifold,
$\nu$ is a smooth immersion,
and the restriction of the symplectic form
$\omega_\Delta$ vanishes on the image
$(d\nu)T_pL$ for any $p\in L$. A topological
map is {\em proper} if the inverse image of
any compact set is compact.

\begin{defn}\label{dLr}
We say that $C\subset\Delta$
is {\em Lagrangian-realizable}
if there exists
a family of proper Lagrangian immersions
\begin{equation}\label{Limm}
\nu_\epsilon:L\to M_\Delta
\end{equation}
smoothly dependent on an arbitrary
small parameter $\epsilon>0$
with the following properties.
\begin{itemize}
\item[(i)] We have
\[
\mu_\Delta(\nu_\epsilon(L))
\subset
C\cup U_\epsilon(V_C).
\]
Furthermore, for each $x\in C\setminus (U_\epsilon(V_C)\cup\dd C\cup\Sigma(C))$
we have 
\[
L\cap \mu^{-1}_\Delta(x)=\{x\}\times (y+\Theta_e)
\subset \mu^{-1}_\Delta(x)=\{x\}\times \Theta,
\]
where $e\subset C$ is the edge containing the point $x$. 
In other words, the intersection 
$L\cap \mu^{-1}_\Delta(x)$ is an affine subtorus
in the fiber $\mu^{-1}_\Delta(x)$.
\item[(ii)] For every $v\in V_\Gamma$ the inverse
image $(\mu_\Delta\circ\nu_\epsilon)^{-1}
(U_\epsilon(v))$
is homeomorphic to the product
$P\times (S^1)^{n-2}$. 
In addition
we have a diffeomorphism of pairs
\[
((\mu_\Delta)^{-1}(U_\epsilon(v)),
\nu_\epsilon(L)
\cap(\mu_\Delta)^{-1}(U_\epsilon(v)))\approx\\
((\C^\times)^2,\phi_v(P_{\delta(v)}))\times (S^1)^{n-2},
\]
where $\phi_v:P_{\delta(v)}\to (\C^\times)^2$
is an embedding whose image is
an irreducible immersed rational holomorphic curve
with three punctures and $\delta(v)$
ordinary nodes.
In particular, all nodes of $\phi_v(P_\delta)
\subset (\C^\times)^2$
are positive self-intersection nodes.
\end{itemize}
\end{defn}

\begin{theorem}\label{main}
Any even primitive tropical curve $C$ in a Delzant
polyhedral domain $\Delta$ is Lagrangian-realizable.
\end{theorem}
This theorem is proved in section \ref{s-proof}.
In the rest of this section we describe
topology of the approximating
Lagrangians $L$ assuming that they exist.
It turns out that their topology is determined
by the tropical curve $C\subset\Delta$ they approximate.

\begin{rmk}
Theorem \ref{main} is expected to be generalized
to Lagrangian realizability of more general
tropical subvarieties (not necessarily curves)
in more general tropical varieties (not necessarily
toric).
In the process of writing the paper I have learned
of a result by Diego Matessi \cite{Ma} establishing
Lagrangian realizability of tropical hypersurfaces
in $\R^n$, $n\le 3$.
In particular, Matessi introduces
the notion of Lagrangian pairs-of-pants for
tropical hypersurfaces in higher dimensions,
which proves to be a very useful new geometric
notion.
Matessi's theorem \cite{Ma} and Theorem 1
share a common special case establishing Lagrangian
realizability for tropical curves in $\R^2$.
\end{rmk}

\ignore{
Let $A_f\subset A$ be a $k$-dimensional affine subspace
of $A$ with an integer slope, i.e. such that
its tangent space $TA_f\subset TA=\Lambda\otimes \R$
is generated by
a $k$-dimensional sublattice of $\Lambda$.
This is equivalent to requiring that the
{\em conormal space} $N^*A_f\subset T^*A=
\Lambda^*\otimes\R$ is generated by an $(n-k)$-dimensional
lattice $\Lambda^*$.

The fiber torus $\Theta=\Lambda^*\otimes\R\approx (S^1)^n$
of the moment map $\mu_\Delta$
can be considered as a (commutative and compact) 
Lie group.
Denote with $\Theta_f\subset\Theta$
the $(n-k)$-dimensional subtorus of $\Theta$ obtained
as $N^*A_f/\Lambda^*$.
If $y\in \Int(M_\Delta)=M_\Delta\setminus
\mu^{-1}_\Delta(\dd\Delta)$ then $y+\Theta_f$
is an affine subtorus of the (torus) fiber of 
$\mu_\Delta:M_\Delta\to\Delta$ containing $y$.
Here we are using the sum as we have
the action of the Abelian group $\Theta$
on $\Int M_\Delta$ coming from 
the action of $\Lambda^*\otimes\R$ on $T^*A$.
}

\begin{lem}\label{Lfeuille}
Let $L\subset M_\Delta$ be a Lagrangian 
subvariety such that 
\begin{equation}\label{cU}
\mu_\Delta(L)\cap U\subset A_f
\end{equation}
for a $k$-dimensional affine subspace $A_f\subset A$
and an open set
$U\subset \Int \Delta$.
Then for any $x\in U\cap A_f$ and
$y\in (\mu_\Delta)^{-1}(x)\cap L$
we have
\begin{equation}\label{yTh}
(\mu_\Delta)^{-1}(x)\cap L\supset
y+\Theta_{f}.
\end{equation}
\end{lem}
\begin{proof}
Since $\Theta_f$ is tangent to the conormal direction
of $A_f$, any of its tangent vector belongs to the
radical direction of the form $\omega_\Delta$ restricted
to $(\mu_\Delta)^{-1}(A_f)$. Since a Lagrangian
subspace is a maximal isotropic direction 
in a tangent space to a symplectic manifold,
any vector parallel to $\Theta_f$ must be contained
in $T_y L$. Thus $(y+\Theta_{f})\cap L$
is of codimension $0$ in $y+\Theta_{f}$ 
for generic $y$
which implies \eqref{yTh}
for all $y\in (\mu_\Delta)^{-1}(x)\cap L$.
\end{proof}

Let $E^b_C\subset E_C$ be the 
set of edges of $C$ of finite
length 
and $e\in E^b_C$.
Choose a point $\iota(e)\in e$ in the
relative interior so that it is disjoint from the (finite)
self-intersection locus of $C$.
Denote $\tilde\mu_\epsilon=\mdn$,
 $T_e=(\tme)^{-1}(\iota(e))$
and $T=\bigcup\limits_{e\in E^b_C}T_e$.

For $v\in V_C\cup\dd C$ we denote by $Q_v$
the component of
$L\sm T$ such that $v\in\tme(Q_v)$
and by $\bar Q_v$ its closure in $L$.
Denote with $\bar P$ a compact pairs-of-pants,
i.e. the complement of three disjoint open
disks in $S^2$. 

For the following series of propositions
we assume that
$C\subset\Delta$ is an even primitive curve,
$\epsilon>0$ is small, and
$\nu_\epsilon:L\subset M_\Delta$
is a Lagrangian immersion satisfying
to the conditions $(i)$ and $(ii)$ of 
Definition \ref{dLr}.

\ignore{
\begin{prop}\label{Ledge}
We have 
\[
T_e=\{\iota(e)\}\times\Theta_e.
\]
\end{prop}
\begin{proof}
The proposition is a direct corollary
of Lemma \ref{Lfeuille}.
\end{proof}
}
Consider the closure $\bar Q_v$ of $Q_v$
in $L$.
By definition of $Q_v$
we have $\dd \bar Q_v=\bigcup\limits_e T_e$,
the union is taken over all $e\in E_C$ adjacent to $v$.
Denote with $\bar P$ the compactification
of the pair-of-pant $P$ into the complement
of three disjoint open disks in $S^2$.
Denote with $\bar P_v\subset\bar P$
a partial compactification of $P$ where we
add a component $E\approx S^1$
of the boundary $\dd\bar P_v$
to $P$ for each {\em bounded} edge $e\in E_C^b$
adjacent to $v$. 
\begin{prop}\label{k0}
A choice of an
$(n-2)$-dimensional affine subspace $A_f\subset A$
(defined over $\Z$)
transversal to $A_v$, $v\in V_C$,
yields a diffeomorphism
\[
\Phi_{v,f}:\bar Q_v\stackrel{\approx}\to
\bar P_v\times\Theta_v
\]
such that for any $x\in P$
we have
\[
\nu_\epsilon(\Phi_{v,f}^{-1}(\{x\}\times\Theta_v))=
y+\Theta_v
\]
for some $y\in M_{\Int \Delta}$,
and for any $a\in \Theta_v$ and
a bounded edge $e\in E_C^b$
adjacent to $v$
there exists $b\in T_e$ with 
\[
\Phi_{v,f}^{-1}(E\times\{a\})=
b+(\Theta_f\cap\Theta_e).
\]
Here $E\subset\dd\bar P$ is the component
corresponding to $e$.
\ignore{
corresponds to .

the $(n-1)$-subtorus $\Theta_e\subset\Theta=T_e$
to $E\times\dd\Theta_v$ for the component
$E\subset\dd\bar P$ corresponding to $e$ so that
the image of
the affine $(n-2)$-subtorus $y+\Theta_v\subset\Theta_e$
for any $y\in \Theta$
is $\{\alpha\times\Theta_v\}$
for some $\alpha\in\dd\bar P$, while the image of
the affine 1-subtorus
$y+\Theta_f\cap\Theta_e\subset\Theta_e$
is $E\times\{\beta\}$ for some $\beta\in\Theta_v$.
}
\end{prop}
\begin{proof}
By Lemma \ref{Lfeuille} the family $y+\Theta_v$
fibers the $n$-dimensional
manifold $\bar Q_v\subset L$
into $(n-2)$-dimensional tori.
Since these tori are affine subtori of
$\Theta=(S^1)^n$, this fibration is 
trivial. Thus its base must be an orientable surface
with boundary $\dd \bar P_v$,
i.e.
a partially compactified pair-of-pants,
perhaps with some handles attached.
By $(ii)$ of Definition \ref{dLr} the base 
must be $\bar P_v$ itself.
Consider a section $\Pi_0$ of the 
(trivial) fibration $\bar Q_v\to\bar P_v$.
Deforming $\Pi_0$ is needed
we may assume that
$\Pi_0(E)$
is an affine subtorus of $T_e$ for each 
$e\in E_C^b$ adjacent to $v$.
The embedding
\[
\Pi_0\subset \bar Q_v\subset M_\Delta\sm
\mu_{\Delta}^{-1}(\dd\Delta)
\]
induces
an embedding of homology groups
\[
\Z^2\approx H_1(\Pi_0)\to H_1(M_\Delta\sm
\mu_{\Delta}^{-1}(\dd\Delta))=\Lambda^*.
\]
The annihilator of its image is
a rank $n-2$ sublattice of $\Lambda$.
Let $A_{f_0}\subset A$ be an $(n-2)$-dimensional
affine subspace parallel to this sublattice.
Thus the trivialization
$\bar Q_v\approx \bar P_v\times\Theta_v$
of the bundle $\bar Q_v\to \bar P_v$
is a diffeomorphism $\Phi_{v,f_0}$
required by the proposition.

Any other $(n-2)$-dimensional direction $A_f$
transversal to $A_v$ corresponds to
the annihilator of
another rank 2 sublattice of $\Lambda^*$
that can be obtained as a graph of a map
$H_1(\Pi_0)\to H_1(\Theta_v)$.
Such a map corresponds to an element of
$H^1(\Pi_0;H_1(\Theta_v))$,
and therefore to a map $\phi:\Pi_0\to\Theta_v$.
As $\Theta_v\subset\Theta$ is a subgroup,
we may use $\phi$ to obtain a new section
\[
\Pi=\{u+\phi(u)\ |\ u\in\Pi_0\}\subset\bar Q_v.
\]
The section $\Pi$ produces a new trivialization
of the bundle $\bar Q_v\to\bar P_v$ and thus
a diffeomorphism $\Phi_{v,f}$ as required by
the proposition. 
\end{proof}
Consider a point $w\in \dd C$.
Let $A_w\subset A$
be the $(n-k)$-dimensional affine
span of the face of $\dd\Delta$
containing $w$. Here $k$ is the codimension
of the boundary point $w$.
Clearly, $\dd\bar Q_w= T_e$, where 
$e\in E_C$ is the edge adjacent to $w$.
\begin{prop}
If $k=2$, and $w$ is a bissectrice point of $\dd C$,
then
the closure $\bar Q_w$, $w\in\dd C$,
is diffeomorphic to $D^2\times (S^1)^{n-2}$.
Under this diffeomorphism $\dd D^2\times\{a\}$,
$a\in (S^1)^{n-2}$
is mapped to the affine 1-subtorus
$y+(\Theta_w\cap\Theta_e)\subset T_e$
for some $y\in T_e$.
\end{prop}
\begin{proof}
By the construction of $M_\Delta$ the inverse image
of the interval $[w,\iota(e)]\subset\Delta$
under $\mu^{-1}_\Delta$ is obtained by taking 
the quotient of $[w,\iota(e)]\times\Theta$
by the 2-torus $\Theta_w$. By Lemma \ref{Lfeuille}
$\bar Q_w\setminus\mu^{-1}_\Delta(w)$
fibers over $(w,\iota(e)]$ with the
fiber $\Theta_e\approx (S^1)^{n-1}$.
The manifold $\bar Q_w$ is obtained by 
contraction of 
the limiting fiber at $w$
by the action of
$\Theta_w\cap\Theta_e\approx S^1$. 
\end{proof}
Suppose now that $k=1$, and the boundary momentum
of $w\in\dd C$ is 2.
This means that the subgroup
$H_w=\Theta_e\cap\Theta_w\subset\Theta$
consists of two elements. Denote the
non-zero element of $H_w$ with $h_w$. 
\begin{prop}\label{k1}
If $k=1$ and the boundary momentum of $C$ at $w\in\dd C$
is 2 then
the closure $\bar Q_w$, $w\in\dd C$,
is diffeomorphic to the (non-orientable)
manifold obtained from $[0,1]\times (S^1)^{n-1}$
by
taking the quotient by the equivalence
$\{0\}\times\{a\}\sim\{0\}\times\{a+h_w\},$
$a\in (S^1)^{n-1}$. In particular,
$\bar Q_w$ fibers over
the M\"obius band with the fiber $(S^1)^{n-2}$.
\end{prop}
Since this equivalence identifies pairs of points
on the boundary of $[0,1]\times (S^1)^{n-1}$
the resulting quotient is a manifold. Since adding
an element of $H_w$ preserves an orientation
of $\{0\}\times (S^1)^{n-1}$, the resulting 
manifold is non-orientable.
\begin{proof}
By Lemma \ref{Lfeuille}
$\bar Q_w\setminus\mu^{-1}_\Delta(w)$
fibers over $(w,\iota(e)]$ with the fiber $\Theta_e$. 
The manifold $\bar Q_w$ is obtained by 
taking the quotient of 
the limiting fiber over $w$
by the action of
$H_w\approx\Z_2$. 

To see that $\bar Q_w$ fibers over the M\"obius
band we choose a $1$-dimensional
subtorus $\Theta_1\subset\Theta_e$
containing the subgroup $H_w$
and a transversal
$(n-2)$-subtorus $\Theta_{n-2}\subset\Theta_e$
such that $\Theta_1\cap\Theta_{n-2}=\{0\}$.
The manifold $\bar Q_w$ fibers over
the M\"obius band obtained from 
$[w,\iota(e)]\times\Theta_1$ by taking
the quotient by the antipodal involution
on $\{w\}\times\Theta_1$.
The fiber is $\Theta_{n-2}\approx (S^1)^{n-2}$.
\end{proof}

\begin{coro}\label{corotop}
We have a diffeomorphism
\begin{equation}\label{topL}
L\approx 
(\bigsqcup\limits_{v\in V_c\cup\dd C}
\bar Q_v)/\sim,
\end{equation}
where the right-hand side is obtained
by gluing the boundaries of the disjoint
union of $\bar Q_v$ according to 
the diffeomorphisms
identifying the components
of $\dd\bar Q_v$ and $T_e$, $e\in E^b_C$,
described by Propositions \ref{k0}-\ref{k1}.
\end{coro}

\begin{rmk}
Note that by Propositions \ref{k0}-\ref{k1},
each component of $L\sm T$
fibers by $(S^1)^{n-2}$. This structure
may be seen as a higher-dimensional generalization
of the {\em graph structure} on 3-manifolds
studied by Waldhausen \cite{Wald}.
In particular, all Lagrangian varieties
produced by Theorem \ref{main} in the case $n=3$
are graph-manifolds. 
\end{rmk}

%% file: exa2.tex
\section{Two-dimensional examples}
\subsection{Lagrangian realizability in the case of planar tropical curves}
Let $\Delta\subset A$ be a Delzant
polyhedral domain in the two-dimensional
tropical affine space $A\approx\R^2$.
Let $C\subset\Delta$ be an even primitive
tropical curve in $\Delta$.
Denote with $j$ the number of boundary points
of $C$ of codimension 1, and by $\kappa$
the number of unbounded edges of $C$.
\begin{thm}\label{th2}
The curve $C$ is Lagrangian-realizable
by a family of Lagrangian immersions
$\nu_\epsilon:L\to M_\Delta$,
for small $\epsilon>0$,
where
$L$ is a connected smooth surface
with $\kappa$ punctures.

If $j=0$ then $L$ is an orientable
surface of genus $b_1(C)$.
If $j>0$ then $L$
is a non-orientable surface homeomorphic to
the connected sum 
of $j+2b_1(C)$ copies of $\rp^2$ with
$\kappa$ punctures.
\end{thm}
This theorem is a special case of 
Theorem \ref{main} and Corollary \ref{corotop}.

We may slightly generalize this theorem to tropical curves
that are not necessarily primitive in $\Delta$ by relaxing
the conditions $(i)$ and $(ii)$ of Definition \ref{pti}.
Namely, let $h_\Delta:\Gamma_\Delta\to\Delta$ be a map
obtained by restriction of a tropical immersion
$h:\Gamma\to\R^2$ to $\Gamma_\Delta=h^{-1}(\Delta)\subset\Gamma$
for a polyhedral domain $\Delta\subset\R^2$.
Assume that $\Gamma_\Delta$ is connected,
$\Sigma(h_\Delta)=
\{x\in\Delta\ |\ \#(h_\Delta^{-1}(x))\ge 2\}$
is finite, and that 
the set $h^{-1}(\dd \Delta)$ is disjoint
from $V_\Gamma\cup\Sigma(h_\Delta)$.
However, instead of $(ii)$ of  Definition \ref{pti}
we only require that for any $e\in E_\Gamma$ the image $dh(e)\in\Z^2$
is non-zero. The greatest common divisor of the coordinates of $dh(e)$
is called the {\em weight} of the edge $h(e)\subset h(\Gamma)$.
Assume that each point
$x\in C\cap\dd\Delta$ sits on an edge of $\Gamma$ of weight 1, and
is either a boundary
point of codimension 1 with the boundary momentum 2,
or a bissectrice point (of codimension 2).

As a planar tropical curve, the curve
$h(\Gamma)\subset\R^2$ is dual to a lattice subdivision ${\mathcal S}_h$
of a lattice convex polygon $N_h$, called {\em the Newton polygon
of $h(\Gamma)$}, see \cite{Mi05}.
Each point $v\in V_\Gamma\cap\Sigma(h)$
may be viewed as a vertex of the rectilinear graph $h(\Gamma)\subset\R^2$,
and corresponds to a subpolygon $N_v\subset N_h$ from ${\mathcal S}_h$.
If $v\in V_\Gamma$ we define
$\delta(v)$ to be the number of lattice points
inside $N_v$. A component of $h(e)\setminus\Sigma(h)$ 
is dual to an edge from ${\mathcal S}_h$. The weight $w(e)$ of 
$e$ is one less than the number of lattice points in the dual edge.
Thus the curve $C=h(\Gamma_\Delta)\subset\Delta$ corresponds 
to a part ${\mathcal S}_{h_\Delta}\subset
{\mathcal S}_h$
formed by the subpolygon and edges of ${\mathcal S}_h$ dual to 
vertices of $C$ and edges of $C\setminus\Sigma(h_\Delta)$.
\begin{figure}[h]
\includegraphics[width=90mm]{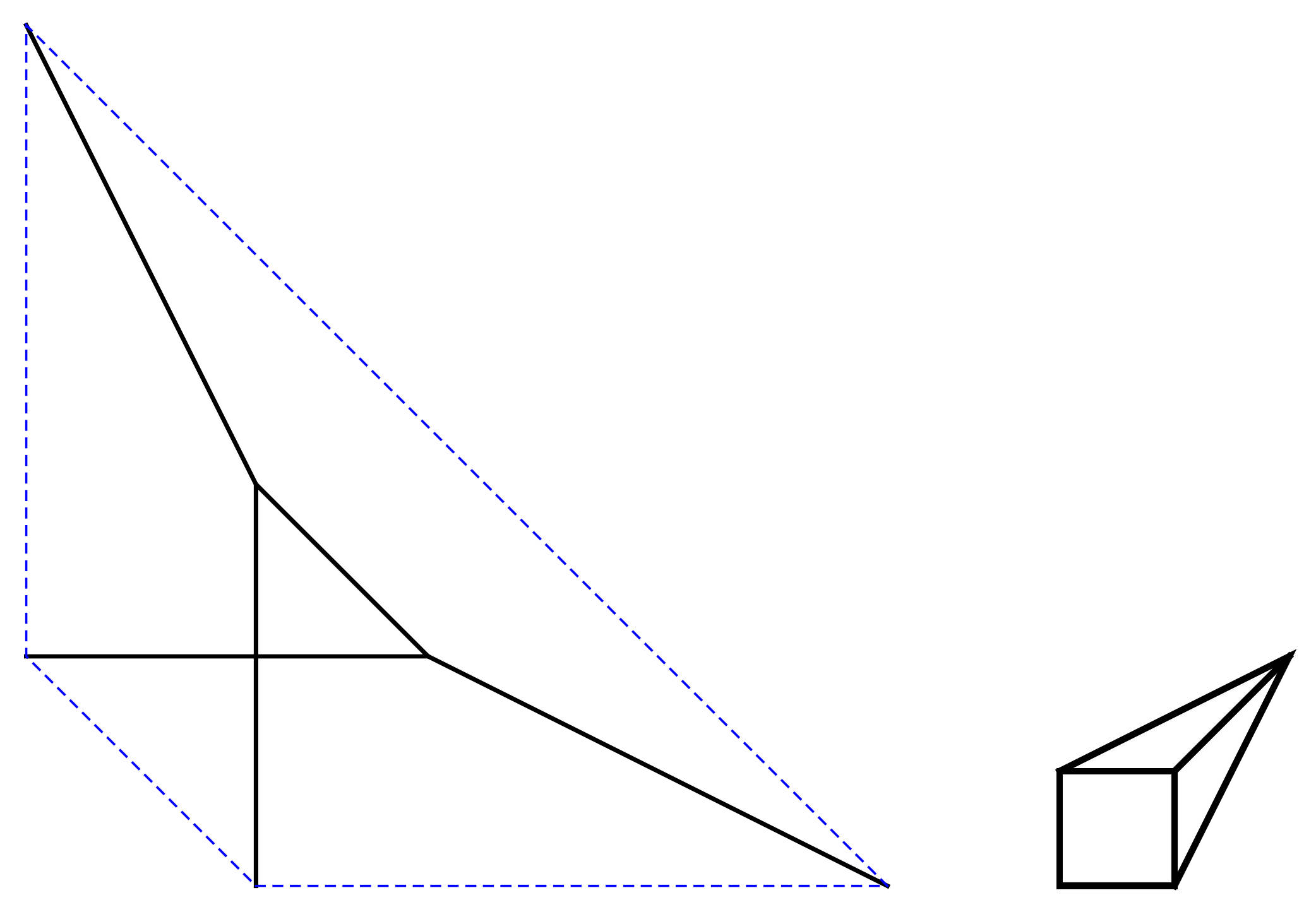}
\caption{An even primitive
tropical curve in a trapezoid
$\Delta$ (dashed), corresponding to an immersed
Lagrangian sphere in $\cp^2\#\overline{\cp^2}$,
and its Newton polygon.\label{fig1}}
\end{figure}

\begin{figure}[h]
\includegraphics[width=90mm]{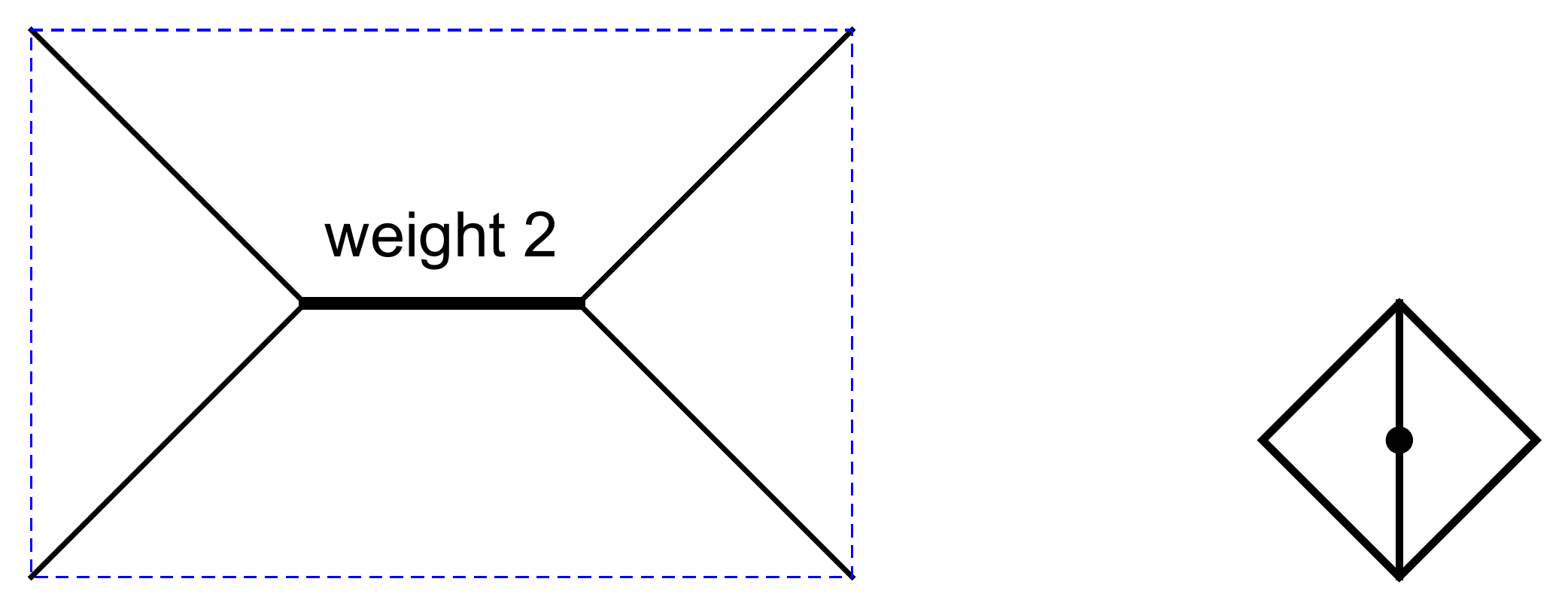}
\caption{A non-primitive
tropical curve in a rectangle for which
Theorem \ref{th2s}
produces an immersed
Lagrangian sphere in $\cp^1\times\cp^1$,
and its Newton polygon.\label{fig2}}
\end{figure}

The curve $C=\Delta\cap h(\Gamma)$ does not have to be
primitive, however we still have
a version of Theorem \ref{th2}. 
As before, we denote the number of codimension 1 boundary
points (of boundary momentum 2) with $j$,
and the number of ends of $\Gamma_\Delta=h^{-1}(\Delta)$ with $\kappa$.
Let 
\[
W=(V_\Gamma\cup\Gamma_\Delta)\cup
\bigcup\limits_{w(e)>1} e
\]
where the latter union is
taken over the edges $e\subset \Gamma_\Delta$
of weight greater than one.
E.g. in Figure \ref{fig2} the set $W$ coincides
with the edge of weight 2.
Note that $W\subset \Gamma_\Delta$
and recall that by our assumption
an edge of $\Gamma$ of weight greater than one must be disjoint
from $\dd\Gamma_\Delta=h^{-1}(\dd\Delta)$.
For simplicity we assume that $V_\Gamma\cap \Gamma_\Delta\neq\emptyset$
(note that we did not make such assumption in Theorem \ref{th2}).

\begin{thm}\label{th2s}
There exists a family of Lagrangian immersions
$\nu_\epsilon:L\to M_\Delta$,
for small $\epsilon>0$,
where $L$ is a connected smooth surface
with $\kappa$ punctures such that 
$\mdn(L)\subset C\cup W_\epsilon$.

Furthermore, let $K\subset W$ be a connected
component and $Y_\epsilon^K$ be the $\epsilon$-neighborhood
of $h(K)$. We have the following properties.
\begin{itemize}
\item
The inverse image $(\mdn)^{-1}(Y_\epsilon^K)\subset L$
contains a unique non-annulus component
$L_{K} \subset L$.
\item 
The surface $L_{K}$ is an (open) orientable surface
of genus $b_1(K)$. The number of ends (punctures) of $L_{K}$
coincides with the number of edges of $\Gamma_\Delta$ adjacent to $K$
plus the number of ends of $K$ itself (if any).
\item
The surface $\nu_\epsilon(L_{K})$ has 
\[
\delta(K)=\sum\limits_{v\in V_\Gamma\cap K}\delta(v)+
\sum\limits_{e\in E_{\Gamma\cap K}}(w(e)-1)+
\#_K(\Sigma(h_\Delta))
\]
ordinary $(+1)$-nodes (i.e. transverse double self-intersection points
with positive local intersection number with respect to 
any orientation of $L_{K}$).
Here $\#_K(\Sigma(h_\Delta))$ is the
number of pairs $x\neq y\in K$ such that
$h(x)=h(y)$, and each such pair is taken
with the weight equal to 
the tropical intersection of the edges $e_x\ni x$
and $e_y \ni y$, i.e. the absolute value
of the scalar product
$(dh(e_x),dh(e_y))$.
\end{itemize}
A component 
$B\subset L\setminus \bigcup\limits_{K\subset W}
L_{K}$
is an annulus if $\mu_\Delta(\nu_\epsilon(B))
\cap\dd\Delta=\emptyset$,
a M\"obius band if $\mdn(B)$ contains a boundary points
of codimension 1, and a disk if $\mdn(B)$ contains a boundary points
of codimension 2.
\end{thm}
This theorem is proved in section \ref{s-proof}
along with Theorem \ref{main}.

\subsection{Example: a real projective plane
inside the complex projective plane}
\begin{figure}[h]
\includegraphics[width=30mm]{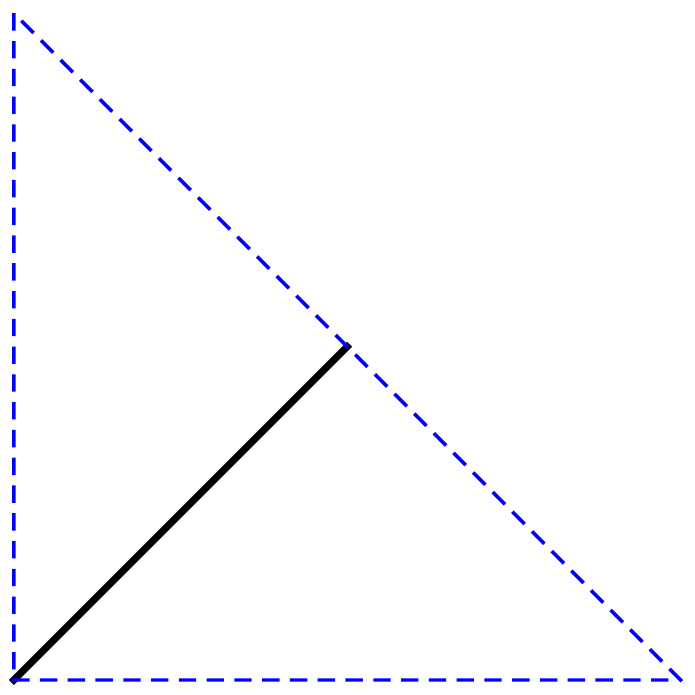}
\caption{One of the ``standard'' Lagrangian
copies of $\rp^2$
in $\cp^2$. \label{fig3}}
\end{figure}

Let $\Delta\subset\R^2$ be the triangle with 
vertices $(0,0)$, $(1,0)$ and $(0,1)$,
and $C$ be the interval between
$(0,0)$ and $(\frac 12,\frac 12)$,
see Figure \ref{fig3}.

The boundary of $C$ consists of two points,
one of codimension 1, and one of codimension 2.
The point $(0,0)$ is a bissectrice point
while the point $(\frac 12,\frac 12)$
is of boundary momentum 2. We have $V_C=\emptyset$.
There is only one edge $e\in E_C$.
In particular, $C$ is smooth.

By Theorem \ref{th2} there exists
a Lagrangian surface diffeomorphic to $\rp^2$
embedded to $M_\Delta=\cp^2$.
Note that since $V_C=\emptyset$ we have
$\mu_\Delta(L)=C$ for any small $\epsilon>0$.

A surface of this kind has appeared
in \cite{CdS-paper} for the
following property.
Consider the interval $I=[(0,a),(1-a,a)]\subset\Delta$
for $0<a<\frac12$.
Then the restriction of $\omega_\Delta$
to $\mu^{-1}_\Delta(I)\approx S^2\times S^1$ is
degenerate,
and has a 1-dimensional radical direction which defines
a $S^1$-fibration $\lambda:\mu^{-1}_\Delta(I)\to M_I=\cp^1\approx S^2$.
We have
$S^1\approx Z=
L\cap \mu^{-1}_\Delta(I)=\{(a,a)\times (y+\Theta_e)$
according to $(i)$ of Definition \ref{dLr}.
The restriction 
$\lambda|_Z:Z\to M_I=\cp^1$
is an embedding as noted in \cite{CdS-paper}.

More generally (by Lemma \ref{Lfeuille}),
we have $L\cap\mu^{-1}_\Delta (t,t)=
\{(t,t)\times (y(t)+\Theta_e)$, where $y(t)\subset\Theta$
may vary with $t\in [0,\frac12]$. If $y(t)=0\in\Theta$ then 
$L$ is the fixed point locus of the antiholomorphic
involution $(z:u:v)\mapsto (\bar z,\bar v, \bar u)$,
and thus is a copy of the standard $\rp^2\subset\cp^2$
in the homogeneous $(z:x:y)$-coordinates 
under the linear substitution $u=x+\bar y$, $v=\bar x + y$.

\subsection{Example: tropical wave fronts of planar polyhedral domains}
Let $\Delta\subset\R^2$ be an arbitrary planar Delzant polyhedral
domain, and $\delta>0$ is small.
The domain $\Delta$ is the intersection of $N$ of half-planes
$\{x\in\R^2\ |\ p_j(x)\ge a_j\}$, $p_j\in\Z^2$, $a_j\in\R$.
$j=1,\dots,N$, where $N$ is the number of sides of $\Delta$.
Without loss of generality we may assume that $p_j$ are
primitive (indivisible) vectors in $\Z^2$.

\begin{figure}[h]
\includegraphics[width=45mm]{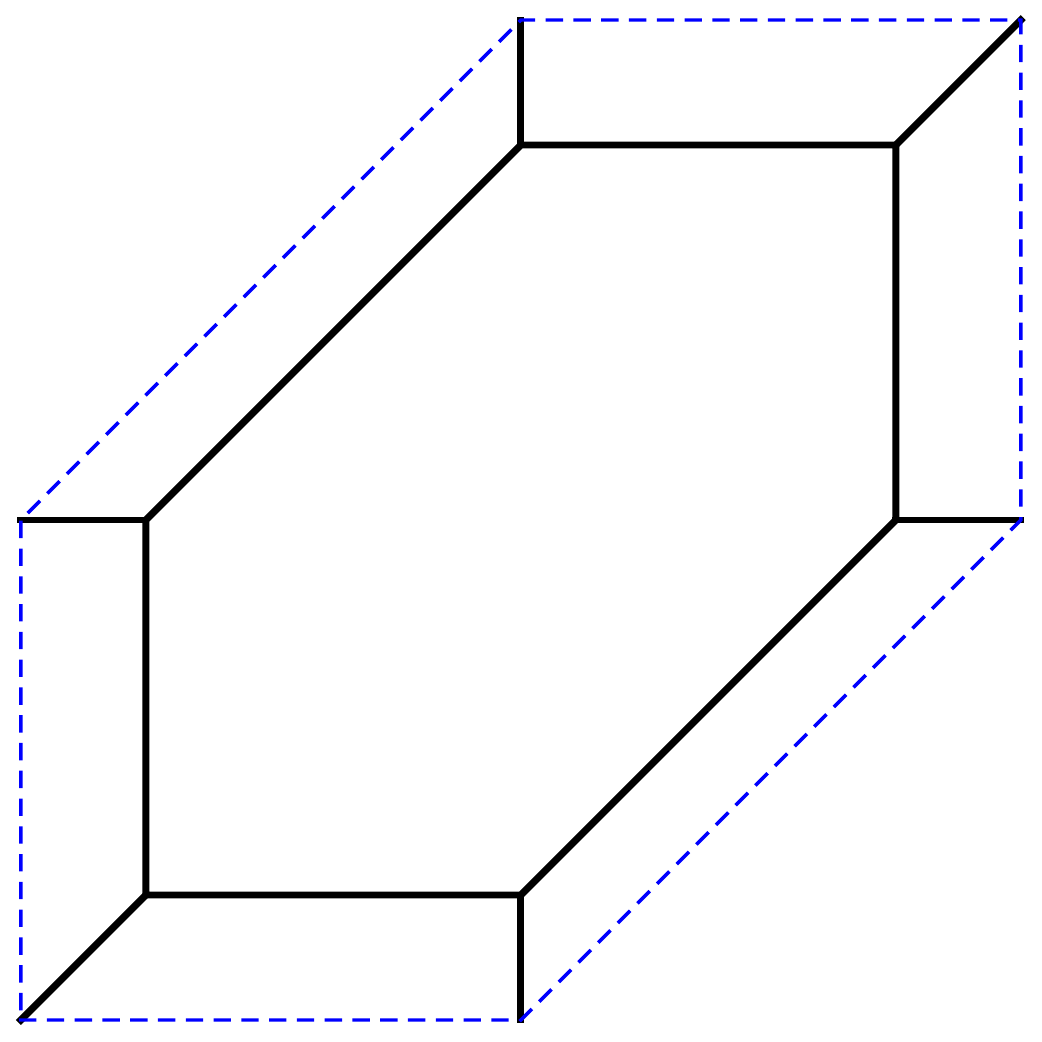}
\caption{Tropical wave fronts of Delzant domains
\cite{KaSh} produce Lagrangian tori in the
corresponding toric surfaces
(the Del Pezzo surface of degree 6 in
the depicted case). \label{fig4}}
\end{figure}

Connecting the vertices of the smaller polyhedral 
domain
\[
\Delta_\delta=\bigcap\limits_{j=1}^N
\{x\in\R^2\ |\ p_j(x) \ge a_j+\delta\}
\]
with the corresponding vertices of $\Delta$,
and taking the union with $\dd\Delta_\Delta$
we get an even primitive tropical curve $W_\delta\subset\Delta$
such that all the vertices of $\Delta$ are the
bissectrice points of $W_\delta$,
see Figure \ref{fig4}.
Since $\Delta$ is assumed to be Delzant, the curve
$W_\delta$ is smooth for small $\delta>0$.

The curve $W_\delta$ can be considered as a tropical
wave front, and appears in the framework related to 
Abelian sandpile models in $\Delta$, see \cite{KaSh}.
In particular, the evolution of this wave front beyond
small values of $\delta$ also produces tropical curves,
though perhaps of different combinatorial type. 

By Theorem \ref{th2} $W_\delta$ is Lagrangian-realizable
by $\Delta$ is Lagrangian-realizable by embedded tori
in the case when $\Delta$ is bounded,
and by embedded cylinders otherwise.

\subsection{Example: Lagrangian embeddings of connected sums of Klein bottles
to $\C^2$}
Recall that
$(\C^2,\frac i2 (dz\wedge d\bar z+dw\wedge d\bar w))$
is the symplectic toric variety corresponding to the
quadrant $\R^2_{\ge 0}$.
Thus even primitive tropical curves in
$\R^2_{\ge 0}$ produce
Lagrangian immersed surfaces in $\C^2$.

\ignore{
The tropical curve
\[
C_3=[(0,0),(1,1)]\cup [(0,3),(1,1)]\cup
[(1,1),(3,0)]\subset\R^2_{\ge 0},
\]
see Figure \ref{c3},
is an even primitive curve with one bissectrice boundary point
and two boundary points of codimension 1 (and of boundary momenta 2).
It is not smooth as
its only vertex has the self-intersection number equal to 2.
Thus we obtain the following proposition. 
\begin{prop}
There exists an immersed Lagrangian Klein bottle in $\C^2$
with two ordinary self-intersection points.
\end{prop}
}

\begin{figure}[h]
\includegraphics[width=90mm]{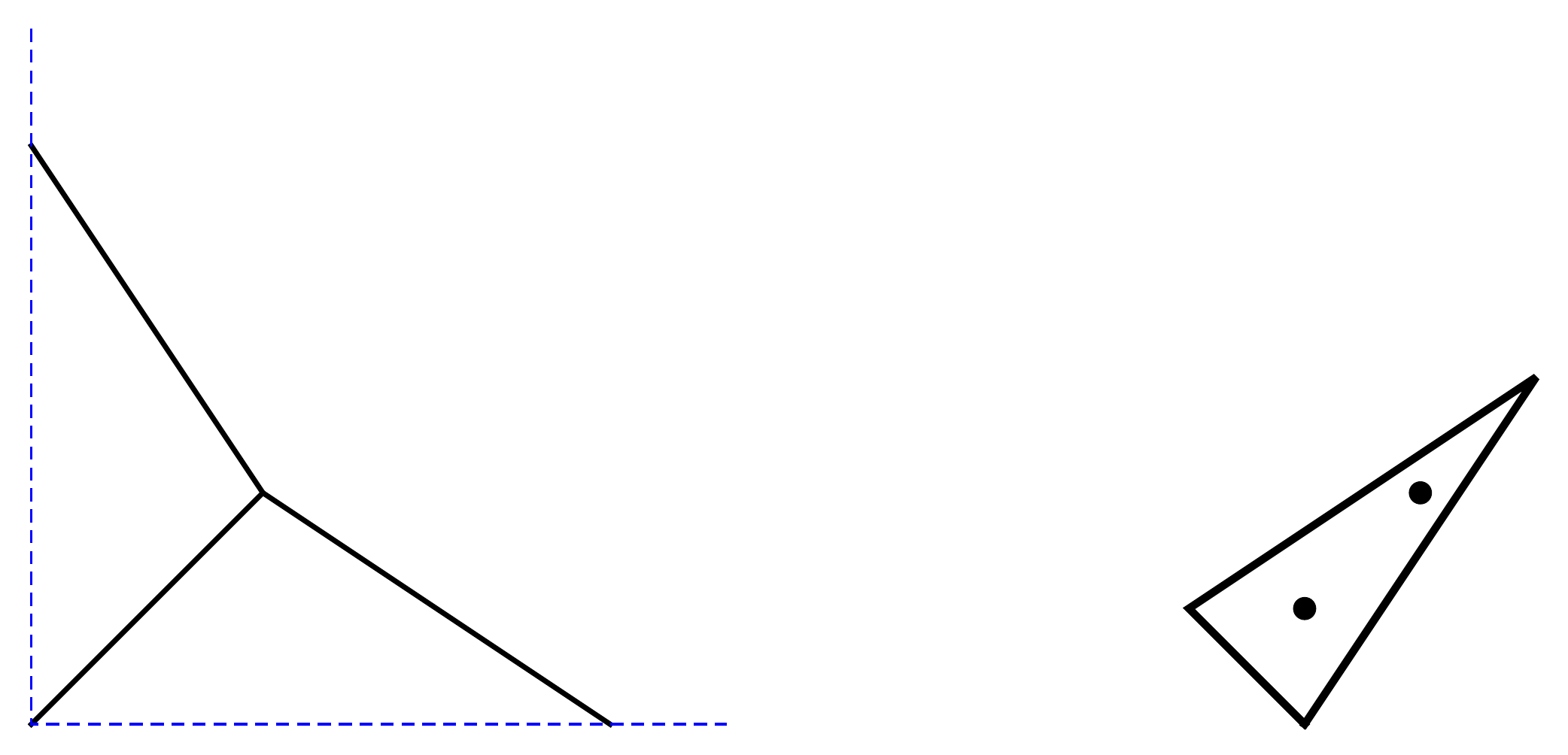}
\caption{A tropical curve representing
an immersed Klein bottle with two nodes
and its Newton polygon. \label{fig5}}
\end{figure}

All boundary points of the tropical curve
\[
C_{3}=[(0,0),(2,2)]\cup [(0,5),(2,2)]\cup [(2,2),(5,0)]\subset\R^2_{\ge 0},
\]
$k>1$,
see Figure \ref{fig5},
are boundary points of codimension 1 and of boundary momenta 2.
The curve is not smooth as
its only vertex has the self-intersection number equal to 2.

Accordingly we get the following proposition.

\begin{prop}\label{kl3}
There exists an immersed Lagrangian Klein bottle in $\C^2$
with two ordinary self-intersection points.
\end{prop}

In a similar way we may construct two series
of Lagrangian immersions of connected sum
of Klein bottles with growing number of self-intersections.

\begin{figure}[h]
\includegraphics[width=90mm]{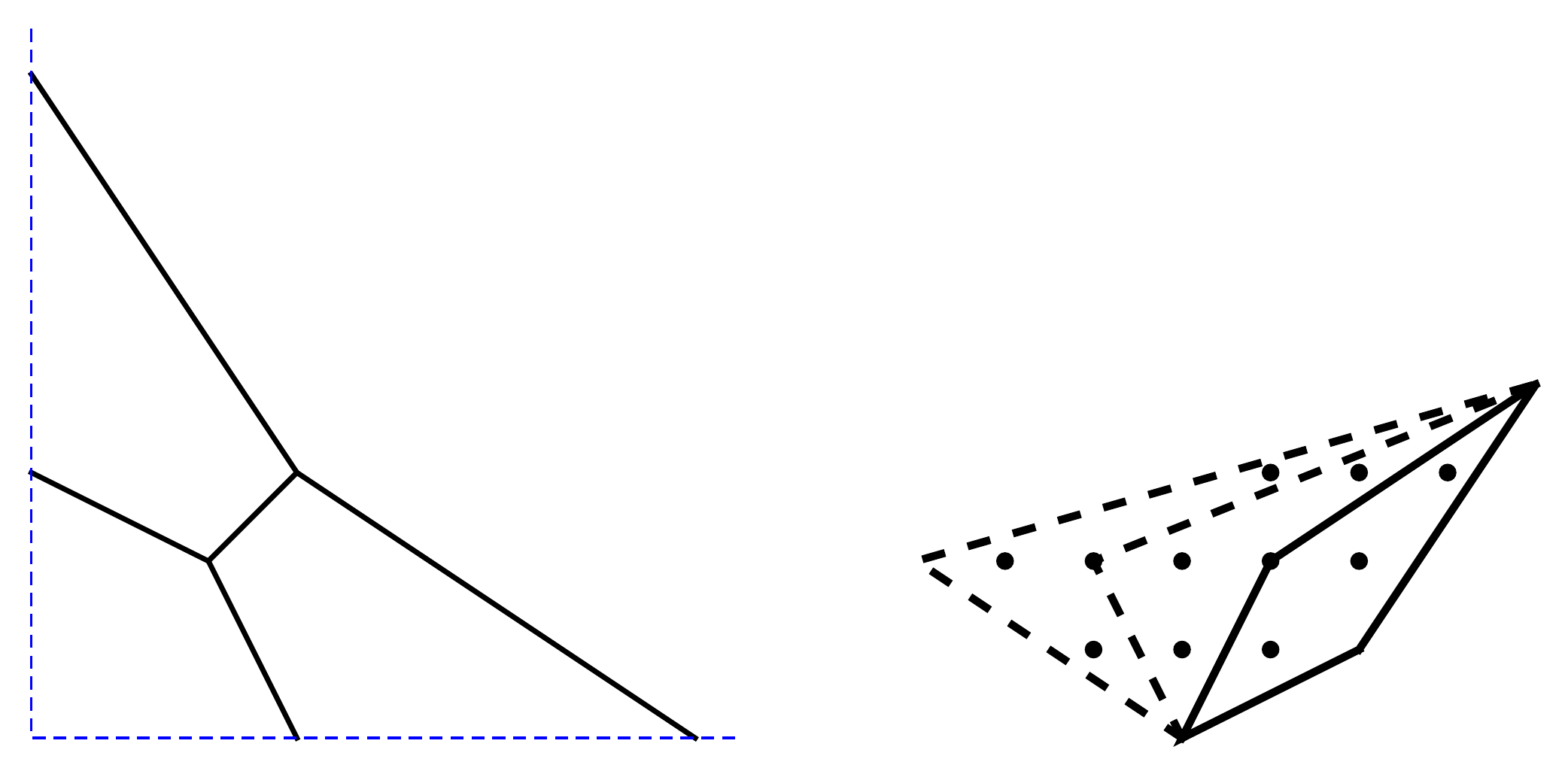}
\caption{The connected sum of two immersed Klein bottle
with $4k-1$ nodes: the tropical curve for $k=1$
and its Newton polygon for $k=1$ (solid polygon) and 
$k=2,3$ (dashed expansion). \label{fig6}}
\end{figure}

\begin{prop}
For any $k\ge 1$ there exist a Lagrangian immersion 
of the connected sum of two Klein bottles to $\C^2$
with $4k-1$ ordinary self-intersection points.
\end{prop} 
\begin{proof}
Consider the tropical curve $C_k$ obtained
as the union of 
the edges $[(0,3),(2,2k)]$,
$[(2,2k),(2+k,0)]$, $[(2,2k),(3,4k-1)]$,
$[(0,7k+\frac 12),(3,4k-1)]$
and $[(3,4k-1),(6k+\frac 32,0)]$.
This is an even primitive tropical curve
in $\R_{\ge 0}^2$
with four boundary vertices of codimension 1,
see Figure \ref{fig6} for $k=1$.
If we extend the boundary edges past $\R_{\ge 0}$ 
as infinite rays we get a tropical curve in $\R^2$
with the Newton polygon
\[
Q_k=\operatorname{ConvexHull}\{(0,-1),(2,2),(-2k+1,0),(-2,-2)\}.
\]
The curve $C_k$ corresponds to the triangulation obtained by
dividing the quadrilateral $Q_k$ into two triangles (corresponding
to the vertices of $C_k$) by the diagonal $[(-2k+1,0),(0,-1)]$, cf. \cite{Mi05}.
The two vertices of $C_k$ have self-intersection numbers $k$ and $3k-1$,
so that Theorem \ref{th2} produces an immersed Lagrangian surface homeomorphic
to the connected sum of Klein bottles with $4k-1$ ordinary self-intersection
points.
\end{proof}
\begin{prop}\label{kl4k}
For any $k\ge 0$ there exist a Lagrangian immersion 
of the connected sum of three Klein bottles to $\C^2$
with $4k$ ordinary self-intersection points.
\end{prop}
\begin{proof}
Consider the lattice triangle $N$ with vertices
$(-2k-1,0)$, $(4,2)$ and $(-4,-2)$.
Consider a subdivision of 
the triangle $N$
into four triangles
by introducing three new vertices
$(-2,-1)$, $(0,0)$ and $(2,1)$ at
the edge $[(4,2),(-4,-2)]$,
see Figure \ref{fig7}.
\begin{figure}[h]
\includegraphics[width=55mm]{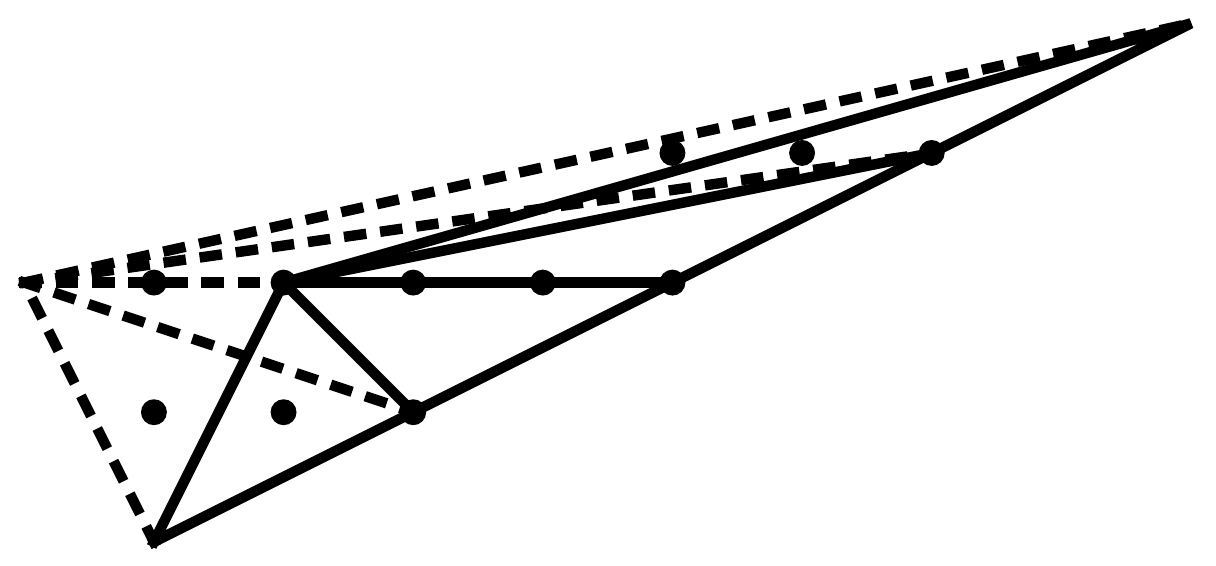}
\caption{Newton polygons and their subdivisions
for Lagrangian immersions 
of the connected sum of three Klein bottles
with $4k$ nodes for $k=1$ and $k=2$. \label{fig7}}
\end{figure}

The resulting triangulation is convex
in the sense of Viro patchworking \cite{Viro-patch}.
Thus there exists a tropical curve $C\subset\R^2$
dual to it, see \cite{Mi05}. By our choice
of the subdivision we have $b_1(C)=0$.
The leaves of $C$ are orthogonal to the edges
of the triangle $N$ serving as the Newton polygon
of $C\subset\R^2$. Namely, we have a leaf in
the direction of $(-2,2k+5)$, a leaf in the direction 
$(-2,3-2k)$, and four leaves in the direction
$(1,-2)$. 

Translating $C$ sufficiently high up in the
direction of $(1,t)$ where $0<t<\frac 2{2k-3}$
ensures that the first two leaves intersect
the $y$-axis of $\R^2_{\ge 0}$ while the last
four leaves intersect the $y$-axis of the translation
so that all boundary momenta are 2. The proposition
now follows from Theorem \ref{th2s}.
\end{proof} 

Smoothing all ordinary nodes of Lagrangian
immersions from Propositions \ref{kl3}-\ref{kl4k}
we recover proof of the following theorem
of Givental.
\begin{coro}[\cite{Gi}]\label{corogi}
There exist Lagrangian embeddings of
surfaces diffeomorphic to $2k+1$ copies
of the Klein bottle to $\C^2$, $k\ge 1$.
\end{coro}

\begin{rmk}
Another construction of Lagrangian embeddings
of $2k+1$ copies
of the Klein bottle to $\C^2$, $k\ge 1$,
can be given by Lagrangian surgery of
an immersed Lagrangian spheres with $2k+1$
self-crossings (out of these $k+1$ must be
positive and $k$ must be negative).
The presence
of negative crossing makes the result
of the surgery non-orientable.
 
By the methods of this paper
it is not clear
if one can further degenerate the surfaces
given by Corollary \ref{corogi} to obtain
such an immersed sphere. 
It might be interesting to detect whether
this deformation exist or Corollary \ref{corogi}
gives a different example of embedded connected
sums of Klein bottles. 
\end{rmk}


%% file: proof.tex
Consider the {\em hyperk\"ahler twist}
\newcommand{\HT}{\operatorname{HT}}
\[
\HT:(\C^\times)^2\to (\C^\times)^2,\
(e^{x_1+iy_1},e^{x_2+iy_2})\mapsto
(e^{x_1-iy_2},e^{x_2+iy_1})
\]
in $(\C^\times)^2$.
\begin{lem}\label{lem-ht}
If $V\subset (\C^\times)^2$ is an immersed holomorphic
curve then $\HT(V)\subset (\C^\times)^2$
is an immersed Lagrangian with respect
to
\[
\omega_{(\C^\times)^2}=\frac i2
(\frac{dz_1}{z_1}\wedge\frac{d\bar z_1}{\bar z_1}+
\frac{dz_2}{z_2}\wedge\frac{d\bar z_2}{\bar z_2})=
dx_1\wedge dy_1+dx_2\wedge dy_2.
\]
\end{lem}
\begin{proof}
The holomorphic 2-form 
$(dx_1+idy_1)\wedge (dx_2+idy_2)$ must vanish
on any tangent space to a holomorphic curve $V$.
Therefore, the form $(dx_1+idy_2)\wedge
(dx_2-idy_1)=(dx_1\wedge dx_2-dy_1\wedge dy_2)
-i(dx_1\wedge dy_1+dx_2\wedge dy_2)$
vanishes on any tangent space to $\HT(V)$.
In particular, the imaginary part on this form
vanishes.
\end{proof}


\ignore{
The following lemma can be considered
as a Lagrangian version of the classical
Menelaus Theorem (see \ref{Mi-AM} for its
tropical version).
}
It is convenient for us to use two
systems of coordinates: the multiplicative
holomorphic coordinates $(z_j)_{j=1}^n
\in (\C^\times)^n$
and the additive real coordinates $(x_j,y_j)_{j=1}^n
\in (\R\times (\R/2\pi\Z))^n$, $z_j=e^{x_j+iy_j}$.
For $t>0$ the scaling map
\newcommand{\sct}{\operatorname{sc}_t}
\begin{equation}\label{scale}
\sct(x_j,y_j)=(tx_j,y_j)
\end{equation}
is a diffeomorphism of $\ctor$ rescaling
the symplectic form $\omega=\sum\limits_{j=1}^n
dx_j\wedge dy_j$. In particular, it sends 
Lagrangian varieties to Lagrangian varieties. 
\begin{lem}\label{lemY}
Theorem \ref{main} holds in the case if $\Delta=\R^2$
and $C\subset\R^2$ is a tropical line,
i.e. $C=Y$ for the union $Y\subset\R^2$
of three rays from $0\in\R^2$
in the direction $(-1,0)$, $(0,-1)$ and $(1,1)$.
\end{lem}
\begin{proof}
The line $V=\{(z_1,z_2)\in (\C^\times)^2\ | \
z_1+z_2-1=0\}$ is holomorphic.
Thus $\HT(V)$ as well as $S_t=\sct(\HT(V))$
are Lagrangians, $t>0$.
In particular, $\mu(S_t)$ is contained
in the $\epsilon$-neighborhood
$U_\epsilon(Y)\supset Y$ for sufficiently
small $t$.

The three rays of $Y$ are symmetric with respect
to $GL_2(\Z)$, so it suffices to deform
$S_t$ into $S'_t$
in the class of Lagrangians within
$\mu^{-1}(U_\epsilon(Y))$ so that 
$\mu(S_t)\cap\{x_1<-\epsilon/2\}=
\{x_2=0\}\cap\{x_1<-\epsilon/2\}$.

But for small $t>0$ over $\{x_1<-\epsilon/2\}$
the Lagrangian surface $S_t$ is approximated by
\[
Z_0=\HT(\{z_2-1=0\})=\HT(\{x_2=0,y_2=0\})=
\{x_2=0,y_1=0\}
\]
since at $z_1=0$ the equation $z_1+z_2-1=0$
degenerates to $z_2=-1=0$.

By the Darboux theorem, a small neighborhood
of $Z_0$ in $\tordva$ is symplectomorphic to 
a small neighborhood in its cotangent bundle
(with the standard symplectic structure $dp\wedge dq$).
Thus any surface approximating $Z_0$ is given
by a small 1-form $\alpha$ on $Z_0$.
As $S_t$ is Lagrangian, we have $d\alpha=0$.
Furthermore, $\alpha$ is exact if and only
if $\int\limits_\gamma\alpha=0$ for a loop
$\gamma$ realizing non-trivial homology class
in $H_1(Z_0)\approx \Z$.

Taking $\gamma=Z_0\cap\{|z_1|=s\}$ for arbitrary
small $s>0$ we see that $|\int\limits_\gamma\alpha|$
must be arbitrary small itself.
Thus $\int\limits_\gamma\alpha=0$,  and so $\alpha=df$
for a smooth function $f:Z_0\to\R$. 
Using the partition of unity
we deform $f$ to a smooth function $f':Z_0\to\R$
such that $f'(x_1,x_2,y_1,y_2)=f$ if $x_1>-3\epsilon/4$,
and  $f'(x_1,x_2,y_1,y_2)=0$ if $x_1<-\epsilon$,
and set $S'_t$ to be the surface defined by $f'$.
\end{proof}

\begin{coro}\label{coroY}
Theorem \ref{main} holds if $\Delta=\R^2$
and $C\subset\R^2$ has a single vertex. 
\end{coro}
\begin{proof}
Since $C$ is trivalent,
there exist an affine map
$\rho:\R^2\to\R^2$ and a multiplicative-affine
map $\rho_{\C}:\tordva\to\tordva$
such that $\mu\circ\rho_{\C}=\rho\circ\mu$
and $C=\rho(Y)$.
Here $\rho$ (resp. $\rho_{\C}$)
is a composition of an
integer linear map in $\R^2$
(resp. in $\tordva$) and a translation in $\R^2$
(resp. in $\tordva$)
of determinant (resp. degree) equal to the multiplicity 
of the vertex $v\in C$, see Lemma 8.21 of \cite{Mi05}.

The image $\rho_{\C}(S_t)$ is a rational curve with
three punctures
of the same Newton polygon as the tropical curve
$C\subset\R^2$, i.e. a pair-of-pants.
Thus it has arithmetic genus $\delta(v)$ 
which correspond to $\delta(v)$ ordinary nodes
by Corollary 8.20 of \cite{Mi05}.
We set $L=\rho_C(S'_t)$ for small
$t>0$ and $S'_t$ from the proof of Lemma
\ref{lemY}.
\end{proof}



\begin{lem}\label{Lconnect}
Let $E_j\in \R^n$, $j=1,2$, be two disjoint 
open intervals
parallel to the same integer vector
$u\in\Z^n$ and
$L_j\subset\ctor$ be two smooth connected
Lagrangian varieties such that
$\mu(L_j)\subset E_j$, and $L_j$ is relatively closed
in $\mu^{-1}(E_j)$.

There exists a Lagrangian
variety $L\subset\ctor$ diffeomorphic to
$\R\times (S^1)^{n-1}$
such that $L_1$ and $L_2$ are its subvarieties 
if and only if $E_1$ and $E_2$ belong to the same line
$E\subset\R^n$.
Furthermore, if $E_1,E_2\subset E$ then the
Lagrangian variety
$L\supset L_1,L_2$ can be chosen so that
$\mu(L)\subset E$.
\end{lem}
\begin{proof}
Without loss of generality we may assume
that $E_j$ are parallel to $(1,0,\dots,0)$.
If $E_1$ and $E_2$ do not belong to the same
line then without loss of generality we may also
assume that the $x_2$-coordinates of $E_1$ and $E_2$
are different. By Lemma \ref{Lfeuille},
both $L_1\cup L_2$ contain $x\times\Theta_E$
for any $x\in E_1\cup E_2$, where
$\Theta_E\subset\Theta=(\R/2\pi)^n$ is the
$(n-1)$-dimensional subtorus conormal to $E$.
 
Let $I=[a_1,a_2]\subset\R^n$ be an interval
connecting $E_1$ and $E_2$.
By Lemma \ref{Lfeuille},
$\mu^{-1}(a_j)=\{a_j\}\times (y_{a_j}+\Theta_E)$.
Choose a smooth function $I\ni a\mapsto y_a\in\Theta=
(\R/2\pi\Z)^n$, and
consider the cylinder
\[
W=\{\{a\}\times (y_{a}+\Theta_2)\ |\ a\in I\}
\subset\ctor.
\]
Here $\Theta_2\subset\Theta$ is a circle
given by $y_k=0$, $k\neq 2$.
Note that the symplectic
area of this cylinder is equal to $2\pi$ times
the difference between the $x_2$-coordinates
of $E_1$ and $E_2$.

Suppose that there exists $L\supset L_1,L_2$
with $H_1(L)\approx\Z$.
Then the two components of $\dd W$ must
be homologous in $L$ and thus the area
of $W$ must vanish.

On the other hand, if $E_1,E_2\subset E$ then
we obtain $L$ by extending
the smooth function $E_1\mapsto E_2\ni a\mapsto
y_a\in\Theta$ to $E$.
\end{proof}

\begin{rmk}
The Lemma above is related to the Menelaus
theorem, in particular, to its tropical version,
see \cite{Mi-AM} for the case $n=2$.
Namely, let $L\subset\tordva$ be a properly
immersed oriented Lagrangian
surface with finitely many ends.
For every end $F$ of $L$ we may take a loop
$\gamma_F$ going around $F$ in the positive direction.
We can define the {\em momentum} $m(F)$ of $F$
to be the integral of the symplectic form 
$\omega$ against a membrane connecting 
$\gamma$ to a loop in the Lagrangian torus
$\mu^{-1}(0)$.
It is easy to see that the momenta of $L$
coincide with the momenta of the corresponding
ends of the tropical curve, and thus
the next proposition corresponds
to Proposition 6.12
of \cite{Mi-AM}. 
\end{rmk}
\begin{prop}[Lagrangian Menelaus theorem]
We have
\[
\sum\limits_F m(F)=0,
\]
where the sum is taken over the ends of
all ends of a properly immersed Lagrangian
surface $L$ in $\tordva$.
\end{prop}
\begin{proof}
The union of all $\gamma_F$ is homologous to zero.
Thus the union of all membranes for $\gamma_f$
can be completed by a Lagrangian chain
(contained in the union of $L$ and $\mu^{-1}(0)$)
to a closed surface in $\tordva$. 
But the integral of $\omega$ against any closed surface
in $\tordva$ is zero.
\end{proof}

\begin{proof}[Proof of Theorem \ref{main}]
Let $C\subset\Delta$ be an even primitive tropical curve
and $v\in V_C$.
The tripod
$Y_v\subset A_v$ obtained as the union
of $v$ with the three rays from $v$ obtained
by extending the edges of $C$ adjacent to $v$
indefinitely is a primitive tropical curve
in the tropical affine plane $A_v\approx\R^2$.
Thus $Y_v\subset A_v$ is Lagrangian-realizable
by Corollary \ref{coroY}.
Let $Q_v$ be the corresponding Lagrangian realization
for $\epsilon/2>0$ restricted to 
$\mu^{-1}(U_\epsilon(v))$.
It is an immersed isotropic
pair-of-pants surface
(which becomes half-dimensional,
and thus Lagrangian once
we restrict the ambient space to a complex monomial
subsurface of $\ctor$
corresponding to $B$).

Consider the $(n-2)$-torus $\Theta_v$ conormal
to all three edges of $C$ adjacent to $v$.
Then $Q_v\times\Theta_v$ gives an immersed Lagrangian
variety $L_v$ in $\mu^{-1}(U_\epsilon(v))$ such
that $\mu(L_v)\cap (U_\epsilon(v)\sm
U_{\frac\epsilon2}(v)$ is contained in $C$.
We use Lemma \ref{Lconnect} for each bounded
edge $[v,v']\in E_C$ to join $L_v$ and $L_{v'}$
to a single Lagrangian.

To finish the proof of the theorem we need
to construct $L$ over the neighborhoods
of boundary points of $C$.
If $w\in\dd\Delta$ is a boundary point
we consider its $\epsilon$-neighborhood
$I_\epsilon(w)\approx [0,\epsilon)$ in $C$
and set $L_w$ to be the closure of
$(I_\epsilon(w)\sm\{w\})\times\Theta_w$,
where $\Theta_w$ is the $(n-1)$-dimensional
subtorus of $\Theta=(\R/2\pi\Z)^n$ conormal
to $I_\epsilon(w)$.
If $w$ is a boundary point of codimension 1
and boundary momentum 1 then $L_w$ is
a smoothly embedded non-orientable Lagrangian variety 
diffeomorphic to the M\"obius band times $(S^1)^{n-2}$.

If $w$ is a bissectrice point
and $\Delta=(\R_\ge 0)^2$
then locally $C$ must coincide with 
the diagonal ray $x_1=x_2$.
Thus the Lagrangian $L_w$ is a smooth 
disk that can be obtained from the complex disk
contained in $\{z_1=z_2\}$ as the image of the
map $\HT$. In the general case $L_w$ is the product
of this disk and $(S^1)^{n-2}$.
\end{proof}

\begin{proof}[Proof of Theorem \ref{th2s}]
Suppose that $\Delta\subset\R^2$, but $C\subset\Delta$
is allowed to have edges of multiple edge 
not adjacent to $\dd\Delta$.
The proof is similar to the proof of Theorem \ref{main},
but instead of patching pairs-of-pants 
over neighborhoods of 
vertices of $C$ we patch together Lagrangian
surfaces over the components of $W$ constructed
in the following way.
Let $K\subset W$ be such a component.
Extend all edges adjacent to $K$ indefinitely
as rays in $\R^2$. The result is a tropical curve.
Its approximation by holomorphic curves 
is provided by Lemma 8.4 of \cite{Mi05} once
we choose a marked point at all but one leaf
of that curve. To finish the proof we apply
the map $\HT$ to the approximating curves.
\end{proof}

%% file: more-bkg.tex
\section{More background material: tropical enumerative
problems and multiplicity}
\subsection{Regular tropical curves}
Given a tropical immersion $h:\Gamma\to A$
we may choose a reference vertex $v\in V_\Gamma$,
and deform $h$ as follows.
We may deform the tropical structure on $\Gamma$,
i.e. vary the lengths of its bounded edges.
Such deformation results in a tropical
curve $\Gamma'$ homeomorphic to $\Gamma$,
but with different lengths of bounded edges.
A tropical immersion $h':\Gamma'\to A$
with $h'(v)=h(v)$ and $dh(e)=dh'(e)$ for
all oriented edges $e\in E_\Gamma$ exists
if and only if 
\begin{equation}\label{eZ}
\sum\limits_{e\in Z}dh(e)l(e)=0\in\Lambda
\end{equation}
for every oriented cycle $Z\subset\Gamma$.

The space of all possible tropical structures
on $\Gamma$ is
$\R_{\ge 0}^b$, where $b$ is the number
of bounded edges of $\Gamma$.
The corresponding deformations of $h:\Gamma\to A$
keeping the value $h(v)\in A$ are 
a subspace $\Def_v(h)\subset\R_{>0}^b$
by the linear equations \eqref{eZ} for a
system of $b_1(\Gamma)$ homologically independent
cycles $Z$.
\begin{defn}
A topological immersion $h:\Gamma\to\R^n$
is called {\em regular} if the codimension of
$\Def_v(h)$ is $nb_1(\Gamma)$, i.e. the conditions
\eqref{eZ} for homologically independent cycles $Z$
are linearly independent.
\end{defn}

Also we may deform $h(v)\in A$ without
changing the tropical structure.
This amounts to composing $h:\Gamma\to A$
with the translation by an element of $\Lambda\otimes\R$
in $A$. If $\Gamma$ is 3-valent then any small
deformation
of $h$ is a combination of a deformation inside
$\Def_v(h)$ and a translation in $A$.
Thus locally the space $\Def(h)$ of all deformations
of $h$ is a linear space
$\Def_v(h)\times(\Lambda\otimes\R)$.
Note that since \eqref{eZ} are linear
equations in $l(e)$ with integer coefficients,
the subspace $\Def_v(h)\subset\R^b_{>0}$
is defined over $\Z$.

From now on we assume that the graph $\Gamma$
with $\kappa$ ends 
is 3-valent, so that $b=\kappa-3+3b_1(\Gamma)$.
The space $\Def(h)=\Def_v(h)\times(\Lambda\otimes\R)$
is a convex open set in a tropical affine space
$\R^b\times\R^n$ of dimension
\[
a=\dim \R^b\times\R^n\ge b-nb_1(\Gamma)+n=\kappa+(n-3)(1-b_1(\Gamma)).
\]
The right-hand side is the expected dimension
of  $\Def(h)$ which agrees with the expected dimension
of the space of deformations of a classical curve
of genus $g$ in an $n$-fold where the value of
the first Chern class on this curve is equal to $\kappa$. 
One may also note
that in any toric compactification of
$\R^n$ the  first Chern class is represented by
the union of divisors at infinity, and so $\kappa$
indeed corresponds to the value of the first Chern
class in the case if the weights of all the leaves
of $h:\Gamma\to\R^n$ are 1.

\begin{rmk}
We say that $h:\Gamma\to\R^n$ is a {\em rational 
tropical curve} in $\R^n$ if $b_1(\Gamma)=0$,
i.e. $\Gamma$ is a tree.
Clearly any such curve is regular as
the codimension in $\R^b\times\R^n$
is tautologically zero.
\end{rmk}

Theorem 1 of \cite{Mi06} claims
that all regular tropical curves are
classically realizable.
The survey \cite{Mi06} did not
contain any proofs, but
various versions of tropical-to-complex correspondences
for regular curves
had appeared
in \cite{Mi05},
\cite{Sh}, \cite{Sp},
\cite{NiSi}, \cite{Tyo}, \cite{Nish}, \cite{MaRu16}
(the list is not exhaustive).
This correspondence is particularly useful
for studying some classical enumerative problems
on the number of curves of given degree and genus
passing through a number of constraints.
For the purposes of this paper we are interested
in the example of a tropical enumerative problem
considered in the next subsection.
In particular, in the remaining part of the paper
we restrict our attention to rational tropical curves.

\subsection{A tropical enumerative problem in $\R^3$
(an example).}
\renewcommand{\DD}{\mathcal D}
\newcommand{\zz}{\mathcal Z}
\newcommand{\qp}{\mathbb Q\mathbb P}
\newcommand{\el}{\mathfrak{l}}
\newcommand{\ez}{\mathfrak{z}}
\renewcommand{\ll}{\mathfrak L}
\newcommand{\ee}{\mathcal E}
Consider the following
example of a tropical enumerative problem in $\R^3$.
Let $\DD=\{d_j\}_{j=1}^\kappa$, $u_j\in\Z^3$,
$\kappa\in\N$, 
be a collection
of non-zero integer vectors
such that $\sum\limits_{j=1}^\kappa d_j=0$.
Such a collection $\DD$ is called a {\em toric degree}.
We say that a (rational) tropical curve $h:\Gamma\to\R^3$
is of degree $\DD$
if  $dh(e_j)=d_j$ for some ordering of the leaves
$\{e_j\}$, $j\in 1,\dots,\kappa$.
We denote the space of all rational tropical
curves of toric degree $\DD$ with $\MM_{\DD}$.

Let $\zz=\{\ez_j\}_{j=1}^\kappa$
$\ez_j\in\Z^3$ be another collection
of non-zero integer vectors, this
time without an extra conditions on the sum.
Consider a configuration $\ll=\{\el_j\}_{j=1}^\kappa$
of lines $\el_j\subset\R^3$ such that $\el_j$ is
parallel to $\ez_j$.
We say that a tropical curve $h:\Gamma\to\R^3$
{\em passes through $\ll$} if
$h^{-1}(\el_j)\neq\emptyset$ for every $j=1,\dots,\kappa$.

Let $\ee=\{e_j\}_{j=1}^\kappa$, $e_j\in E_\Gamma$,
be a collection of edges of $\Gamma$ (maybe
with repetitions).
We say that $h:\Gamma\to\R^3$ {\em passes
through $\ll$ at $\ee$} if 
$h^{-1}(\el_j)\cap e_j\neq\emptyset$\
for every $j=1,\dots,\kappa$.

For the collection $\zz$ we denote the space of
all compatible configurations $\ll$
with $\MM(\zz)$.
Clearly, $\MM(\zz)\approx (\R^2)^\kappa$
is an affine space over a vector space of
dimension $2\kappa$. 
Thus we may speak
of generic configurations $\ll$.

\begin{defn}
We say that $h:\Gamma\to\R^3$ and 
$h':\Gamma'\to\R^3$ are {\em of the same combinatorial
type} if there exists a homeomorphism
$\Phi:\Gamma\to\Gamma'$ that sends
an oriented edge $e\in E_\Gamma$ to 
an oriented edge $e'\in E_{\Gamma'}$
so that $dh(e)=dh'(e')\in\Z^3$.
We have denoted the combinatorial type of $h$
with $\Def(h)$.
\end{defn}

Clearly, there are finitely many combinatorial
types of curves of a given toric degree
passing through $\ll$ at $\ee$.

\begin{lem}
If $\ll$ is chosen generically in $\MM(\zz)$
then
there is at most one curve in the combinatorial
type of a tropical rational curve
of toric degree $\DD$
passing through $\ll$ at $\ee$.

Furthermore, if a tropical curve
$h:\Gamma\to\R^3$ passes through $\ll$
then $\Gamma$
of $h$ is 3-valent, and the inverse image
$h^{-1}(\el_j)$
is finite and disjoint from $V_\Gamma$.
In addition we have $dh(e)\neq 0$ for any
$e\in E_\Gamma$.
\end{lem}
\begin{proof}
By regularity of rational tropical curves we have
$\dim\Def(h)\le\kappa$, and $\dim\Def(h)<\kappa$ 
unless $\Gamma$ is 3-valent.
For a vertex $v\in V_\Gamma$ consider
the map $\ev_v:\Def(h)\to\R^3$ defined
by $\ev_v(h)=h(v)$.
For an edge $e\in E_\Gamma$ define
the space $\R^3/dh(e)\approx\R^2$
as the quotient of $\R^3$ by the subspace
generated by the vector $dh(e)$ if it is
non-zero and as $\R^3/dh(e)=\R^3$ if $dh(e)=0$.
Define the map
\begin{equation}\label{eve}
\ev_e:\Def(h)\to\R^3/dh(e)
\end{equation}
by setting $\ev_e(h')$, $h'\in\Def(h)$,
to be the image of the edge
$h'(e)$ under the projection $\R^3\to\R^3/dh(e)$.
Clearly, these are linear maps on 
$\Def(h)\subset\R_{>0}^b\times\R^3$.
The lemma now follows from computing dimensions
once we note that if $dh(e')=0$ for an edge $e'$
then the map $\Def(h)\to\R^3/dh(e)$ has a kernel
corresponding to varying of the length of $e'$
for any $e\in E_\Gamma$.
\end{proof}

\begin{coro}
There is a finite set $\MM(\DD,\ll)$ of rational tropical curves
of toric degree $\DD$ passing through $\ll$
as long as $\ll$ is generic.
\end{coro}

For every $h:\Gamma\to\R^3$ such that
$h^{-1}(\el_j)\cap e_j\neq 0$ define
\begin{equation}\label{eqeve}
\ev_{\ee}:\Def(h)\to\R^l
\end{equation}
as the direct sum of the maps
sending $h'\in\Def(h)$
to the mixed product of $d_j$,
$\ez_j$ and any vector connecting a point
of $h(e_j)$ and a point of $\el_j$.
The curve $h:\Gamma\to\R^3$ passes through $\ll$ at $\ee$
if and only if $\ev_{\ee}=0$. 
The map $\ev_{\ee}$ is linear and defined over $\Z$.

The {\em multiplicity $m(h,\ee,\ll)$} of a tropical curve $h$
passing through $\ll$ at $\ee$ is the absolute
value of the determinant of $\ev_e$.
We set 
\[
m(h,\ll)=\sum\limits_{\ee}m(h,\ee,\ll),
\] 
where the sum is taken over all possible
configuration of edges $\ee=\{e_j\}$, $e_j\in E_\Gamma$,
$j=1,\dots,k$, with
$h^{-1}(\el_j)\cap e_j\neq\emptyset$.
The number
\begin{equation}\label{ntrl}
m(\DD,\ll)=\sum\limits_{h\in\MM(\DD,\ll)}m(h,\ll)
\end{equation}
is
{\em the number of tropical curves passing through
$\ll$}.

It can be shown that $m(\DD,\ll)$ depends only
on $\DD$ and $\zz$ (recall that $\ll$ 
is chosen generically). In particular,
this claim follows from the correspondence with 
the number of complex curves
in $(\C^\times)^3$
of toric degree $\DD$ passing through an
appropriate collection of monomial curves
defined by $\zz$
which is a special case of \cite{NiSi}.

%% file: l3.tex
\newcommand{\mv}{\operatorname{mv}}
\subsection{A vector product formula for
the tropical multiplicity $m(h,\ll)$}
In this subsection we compute 
the multiplicity $m(h,\ll)$ from
\eqref{ntrl} for
a rational 3-valent tropical curve $h:\Gamma\to\R^3$
of toric degree $\DD$
passing through a generic configuration $\ll$
of lines $\el_j\subset\R^3$ parallel 
to a collection $\zz=\{\ez_j\in\Z^3\}$,
$j=1,\dots,\kappa$ in a special case
of the so-called
{\em boundary configuration}. 
\begin{defn}
We say that $\ll$ is a boundary configuration
for $h$ if for every $j=1,\dots,\kappa$
there exists a leaf (unbounded edge)
$e_j\subset\Gamma$ such that 
$h^{-1}(\el_j)$ is a point contained in
$e_j$ and $e_j\neq e_k$ if $j\neq k$.
\end{defn}
For the rest of the section
we assume that $\ll$ be a boundary configuration
for $h$. 
Recall that $d_j=dh(e_j)\in\Z^3$ is
the image of a unit tangent vector
to the edge $e_j$ under
the differential of $h$.
\begin{defn}
The {\em leaf rotational momentum}
\[
\rho_j=d_j\times\ez_j\in\Z^3
\]
is the vector product of $u_j$ and $\ez_j$.
Like the vectors $u_j$ and $\ez_j$ themselves,
it is well-defined up to sign.
\end{defn}

\begin{figure}[h]
\includegraphics[width=50mm]{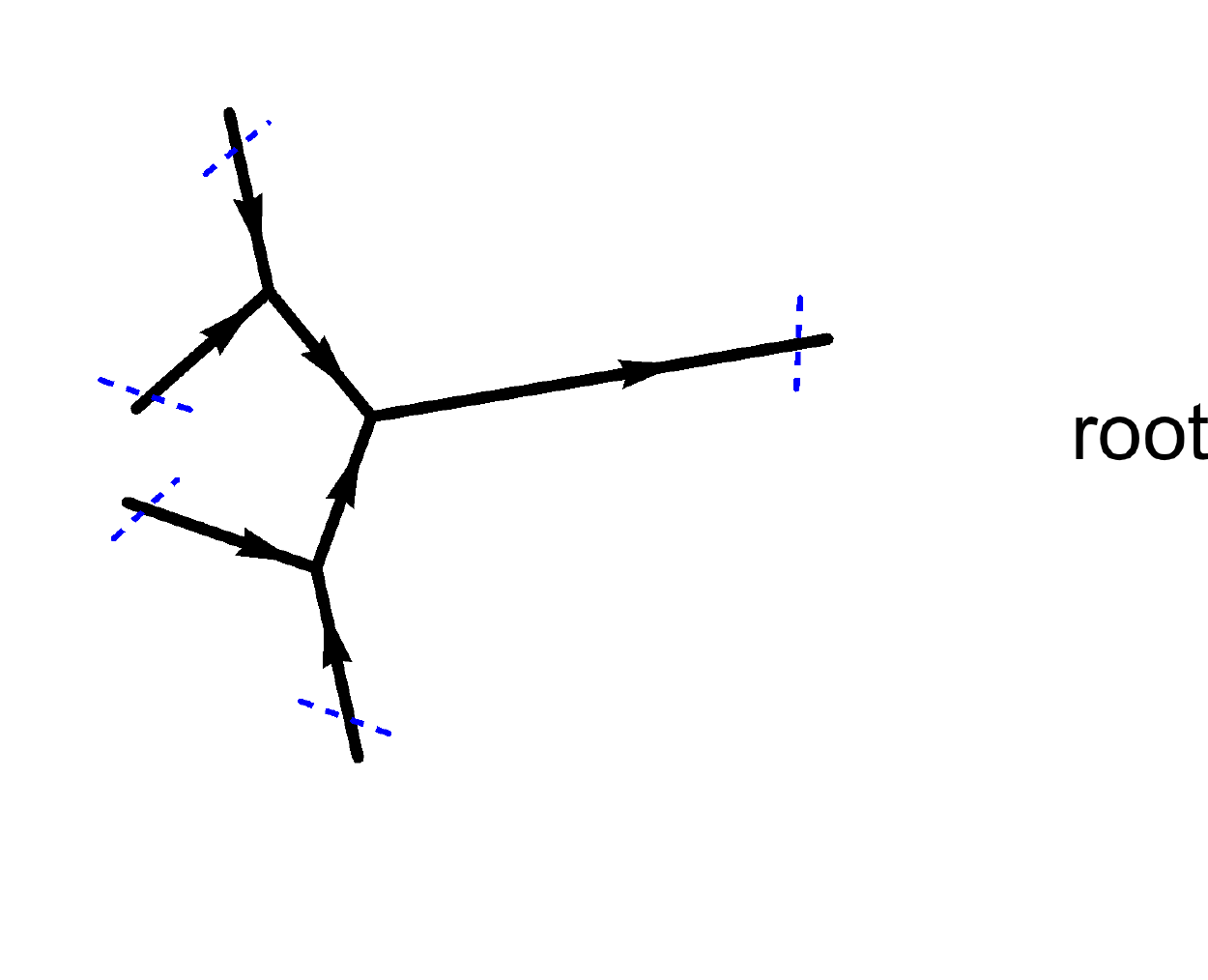}
\caption{A tropical curve trough a boundary
configuration of lines (shown by dashed intervals)
and the order of binary
operations determined by a choice of 
the root. \label{fig8}}
\end{figure}

The tree $\Gamma$ with the leaf $e_\kappa$
chosen as the root
encodes a way to place parentheses for the
binary operation on the leaves $e_1,\dots,
e_{\kappa-1}$.
Namely, $\Gamma$ is oriented towards the root,
and the orientation determines the order
of binary operations, see Figure \ref{fig8}.
We start with the rotational 
momenta at the leaves $e_1,\dots,e_{\kappa-1}$.
At every vertex of $\Gamma$ 
we have two incoming edges $\iota_1,
\iota_2\in E_\Gamma$, and an outgoing edge
$o\in E_\Gamma$.

If the rotational momenta
for the incoming edges are already
defined then we define the rotational momentum for the
outgoing edge as
\begin{equation}\label{rhoo}
\rho(o)=(\rho(\iota_1)\times\rho(\iota_2))
\times (dh(o)),
\end{equation}
i.e. the vector product of the rotational momenta
of the incoming edges followed by the vector
product with the image of a unit vector to the
outgoing edge.
The vector $\rho(o)$ is well-defined
up to sign as we do not have a preferred order
for the incoming edges, 
but so are the 
rotational momenta for the incoming edges. 
\begin{rmk}\label{r-rhoo}
Presence of the triple vector product in 
\eqref{rhoo} is related to our usage of a vector
product which depends on the scalar product
in $\R^3$ that is not canonical for our problem.
However, thanks to the double appearance
of the vector product operation, the result
depends only on the volume form in $\R^3=\Z^3\otimes\R$
that is naturally associated with the tropical
structure, cf. Remark \ref{r-bivp}.
Namely, the image of the unit vector
$dh(o)$ and the rotational momenta $\rho(\iota_j)$
belong to dual vector spaces. More canonically,
we may replace the vector product with the wedge 
product followed by the isomorphism to the
dual space (cf. \eqref{bivp}), so that the image of
$\rho(\iota_1)\wedge\rho(\iota_2)$ and $dh(o)$ 
belong to the same (dual) vector space. Taking
their wedge product followed by the isomorphism
back to the original space we obtain the rotational
momentum $\rho(o)$ from the same vector space as
$\rho(\iota_j)$.
\end{rmk}

Starting from $\rho(e_j)$, $j=1,\dots,\kappa-1$,
and applying \eqref{rhoo} at
all vertices we obtain the rotational momentum
$\rho(e_\kappa)$.
\begin{defn}\label{def-mhp}
The {\em mixed $h$-product}
\[
(\ez_1,\dots,\ez_\kappa)_h = \rho(e_\kappa).\rho_\kappa
\]
of rotational momenta
of $\zz$ along $h:\Gamma\to\R^3$
is the scalar product of $\rho_\kappa$ and
$\rho(e_\kappa)$. It is well-defined up to sign.
\end{defn}

\begin{exa}\label{exaY}
If $\Gamma$ is a tripod then $\kappa=3$
and 
\[
(\ez_1,\ez_2,\ez_3)_h=\pm (\rho_1,\rho_2,\rho_3)
\in\Z^3
\]
is the mixed product
of the rotational momenta $\rho_1,\rho_2,\rho_3\in\Z^3$.
\end{exa}

\begin{lem}
The $h$-product $(\ez_1,\dots,\ez_\kappa)_h$ does not
depend on the choice of the root leaf $e_\kappa$
and thus on the order of $\ez_1,\dots,\ez_\kappa$.
\end{lem}
\begin{proof}
It is convenient to allow a 3-valent vertex
$v\in\Gamma$ to be a root as well.
We have the binary
operations for rotational momenta as before
for edges with two input and one output vertices.
At $v$ we have three input vertices and
perform the mixed product.

By the property of the mixed product 
of three vectors in $\R^3$ the result (up to sign)
is invariant under all permutations
of these vectors.
This implies that the mixed $h$-product
stays invariant if we change the root from
a leaf to its adjacent 3-valent vertex,
or from a 3-valent vertex to an adjacent
3-valent vertex.
\end{proof}

\begin{prop}\label{prop-vprod}
The absolute value $|(\ez_1,\dots,\ez_\kappa)_h|$
coincides with the multiplicity $m(h,\ll)$.
\end{prop}
\begin{proof}
The proposition amounts to computing the determinant of the map
\eqref{eqeve} in the case when $\Gamma$ is a trivalent tree,
and the line $\el_j$ intersects the leaf $e_j$ of $\Gamma$.
We do it by induction on the number of vertices of $\Gamma$.

If $\Gamma$ has a single 3-valent vertex then
$\Def(h)$ can be identified with translations in $\R^3$.
The coordinates of $\ev_{\ee}$ are given by the rotational momentum
of the leaves $e_j$.
Thus the proposition in this case follows from the interpretation
of the (classical) triple mixed product as the determinant of
the parallelepiped built on the vectors.

Let $\Gamma$ be a tree with at least two 3-valent vertices,
and $p\in\Gamma$ be a point inside a bounded edge $e$
of $\Gamma$. 
The complement $\Gamma\sm\{p\}$ consists of two components.
Extending their edges adjacent to $p$ indefinitely to rays $r_j$
we obtain two rational 3-valent tropical curves
$h_j:\Gamma_j\to\R^3$, $j=1,2$, 
with the property that $h(r_1)$ and $h(r_2)$ are contained in
the same line of $\R^3$.
We get a linear inclusion map $\Def(h)\to\Def(h_1)\times\Def(h_2)$
with the image given as the pull-back of the diagonal in
$(\R^3/dh(e))\times (\R^3/dh(e))$
under the map
\[
\ev_{r_1}\times\ev_{r_2}:\Def(h_1)\times\Def(h_2)\to 
(\R^3/dh(e))\times (\R^3/dh(e))=E\times E.
\]
The vector space $E=(\R^3/dh(e)) \approx\R^2$ is by construction
defined over integers. Denote the underlying integer lattice
with $\Lambda_E\approx\Z^2$.

The generator $\Lambda_h\in\Lambda^\kappa(\Def(h))\approx\R$
of the top exterior power of the integer lattice in $\Def(h)$
is well-defined up to sign.
It is given as the pull-back of the integer generator of
the exterior square
of the diagonal in $E\times E$ which in its turn is given as 
\begin{equation}\label{ddelta}
\delta^0_1\otimes\delta^2_2+\delta^2_1\otimes\delta^0_2+
\delta^{a}_1\otimes\delta^b_2+\delta^{b}_1\otimes\delta^a_2,
\end{equation}
where $\delta^i_j\in\Lambda^i(\Lambda_{E_j})$, $i=0,2$, 
and $\delta^a_j,\delta^b_j\in\Lambda^1(\Lambda_{E_j})$
are given by an integer basis $a,b\in\Lambda_{E_j}$,
where $E_1$ (resp. $E_2$) is the first (resp. the second)
copy of $E$ in the product $E\times E=E_1\times E_2$.


Note that the number of vertices of each tree $\Gamma_j$
is smaller than that of $\Gamma$. 
Let us choose $r_1$ as the root leaf
of $\Gamma_1$ and compute the rotational momentum
$\rho(r_1)$ starting from all the other leaves of $\Gamma_1$. 
Let $\el_a$ be the line parallel to
$\rho(r_1)$ and passing through $p$
and $\el_b$ be the line parallel to
$\rho(r_1)\times (dh_1(r_1))$ and passing through $p$.
Note that the projections of $\el_a$ and $\el_b$ to $E$
are two lines containing a basis of $\Lambda_E$.

We claim that
\begin{equation}\label{mprod}
m(h,\ll)=
m(h_1,\ll_1\cup \el_a)m(h_2,\ll_2\cup \el_b)/w(e),
\end{equation}
where $\ll_j$ is the subconfiguration of $\ll$
consisting of lines adjacent to $h_j(\Gamma_j)$,
and $w(e)$ is the weight of the edge $e\in E_\Gamma$,
i.e. the ratio of $dh(e)\in\Z^3$ and a
primitive integer vector in the same direction.

To see this we use \eqref{ddelta} and use a basis of $\Lambda^1(\Lambda_E)$
given by the projection of the lines $\el_a$ and $\el_b$. 
The evaluation map \eqref{eqeve} extends to the ambient space
$\Def(h_1)\times\Def(h_2)$ that corresponds to disconnected
graphs $\Gamma_1\cup\Gamma_2$. The pull-back of the integer volume
form on $\R^\kappa$ vanishes on the first two summands of \eqref{ddelta}
by dimensional reasons.
Furthermore, it vanishes on the last summand of \eqref{ddelta} by
the definition of $\rho(r_1)$ and the induction hypothesis applied to $h_1$
as the scalar product of $\rho(r_1)$ and $\rho(r_1)\times (dh_1(r_1))$
is zero, and thus $m(h_1,\ll_1\cup\el_b)=0$.
The remaining term corresponds to
$m(h_1,\ll_1\cup \el_a)m(h_2,\ll_2\cup \el_b)/w(e)$
as the factor $w(e)$ appears two times in $\ev_{r_1}\times\ev_{r_2}$.
But the projection of $\rho(r_1)$ to $E$ is $\pm m(h_1,\ll_1\cup \el_a)a$
and the proposition follows from the induction hypothesis applied to $h_2$.
\end{proof}

\begin{rmk}\label{r-bivp}
As we have already seen (cf. Remark \ref{r-rhoo}),
the vector product presentation of tropical
multiplicity of Theorem
\ref{m-thm} is based on the explicit
isomorphism given between
$\R^3\approx V=\Lambda\otimes\R$
and its dual vector space
$\R^3\approx V^*=\Lambda^*\otimes\R$
provided by the scalar product. Alternatively,
a vector product in $V$
can be replaced with the wedge
product followed by the isomorphism between
the wedge square of $V$
and $V^*$
provided by the volume form coming from the 
integer lattice $\Lambda\subset V$ (see \eqref{bivp}).
Note that the latter isomorphism is canonical 
up to a sign. This viewpoint may be used
for generalization to higher dimensions.

As the author has learned during the final stage of 
preparation of this paper, the vector product formula
for tropical multiplicities from Proposition \ref{prop-vprod}
is generalized to higher dimensions in the upcoming work of Travis Mandel and
Helge Ruddat \cite{MaRu} to polyvector calculus replacing the vector product.
Their generalization also allows working with curves of positive
genus as well as with gravitational descendants.

Also an earlier work of Mandel and Ruddat
\cite{MaRu16} has
announced appearance of tropical
Lagrangian lens spaces $L(p,q)$ in
the mirror dual quintic 3-folds associated
to tropical lines of multiplicity
$p$ which coincides with $\#(H_1(L(p,q)))$
in an upcoming work of Cheuk Yu Mak and
Helge Ruddat \cite{MakRu}. Their work
should be particularly relevant
to Theorem \ref{m-thm} and Example \ref{exa-ls}.
%
\end{rmk}

\subsection{Lagrangian rational homology
spheres in toric 3-folds}
Let $C$ be a compact even primitive
tropical curve in a Delzant polyhedral
domain $\Delta\subset\R^3$ which does not
have boundary points
of codimension 1 (i.e. such that all of its
boundary points are bissectrice points).
By Theorem \ref{main} the curve $C$
is Lagrangian-realizable
by immersions
$\nu_\epsilon:L\to M_\Delta$
of oriented closed
(compact without boundary)
smooth 3-manifolds $L$.
Suppose that $C$ has $\kappa$ bissectrice
points $w_j$, $j=1,\kappa$.
Each $w_j$ is 
is adjacent to an edge $e_j\subset C$.
Extending edges $e_j$ to unbounded rays 
through the exterior of $\Delta$ gives us
a tropical curve $h:\Gamma\to\R^3$
such that $C=h(\Gamma)\cap\Delta$.
Also $w_j$ is contained in an edge 
of the 1-skeleton of $\dd\Delta$.
Let $\el_j\subset\R^3$ be the line 
containing this edge and $\ll=\{\el_j\}_{j=1}^\kappa$.
Choose $\ez_j\in\Z^3$ to be one of the two
primitive vectors parallel to $\el_j$. 
The the mixed $h$-product $(\ez_1,\dots,\ez_\kappa)_h$
is now given by Definition \ref{def-mhp}.
Define the {\em vertex multiplicity} of $C$ as
$\mv=\prod\limits_{v\in V_C}m(v),$
where $m(v)$ is defined by \eqref{mult-v}.
 
\begin{thm}\label{m-thm}
If $(\ez_1,\dots,\ez_\kappa)_h\neq 0$
then for small $\epsilon>0$ 
the $3$-manifold $L$ is an oriented rational 
homology sphere such that
\[
\#(H_1(L))=
\frac{m(h,\ll)}
{\mv}
=\frac{|(\ez_1,\dots,\ez_\kappa)_h|}
{\mv}.
\]
\end{thm}
Let $\Delta_t$, $t\in [0,1]$ be a small deformation
of the polyhedral domain $\Delta$.
This means that
\[
\Delta(t)=\bigcap_{j=1}^N\{x\in\R^3\ |\
p_jx\ge a_j(t)\}\subset\R^3,
\]
$\Delta(0)=\Delta$,
$p_j\in\Z^3$ and $a_j(t)\in\R$ is a
slowly varying
smooth function.
\begin{coro}\label{m-coro}
If $(\ez_1,\dots,\ez_\kappa)_h\neq 0$
then for small $\epsilon>0$ there exists a 
smooth deformation
\[
\nu_{\epsilon}(t):L(t)\to M_{\Delta(t)},
\]
$t\in [0,1]$.
\end{coro}
More precisely, the symplectic quotient construction
applied to the deformation $\Delta_t$ 
yields a topological space $M_{\tilde\Delta}$
that maps to $[0,1]$ with the fiber $M_{\Delta_t}$
over $t\in [0,1]$.
Note that $M_{\tilde\Delta}$ is a smooth manifold
outside of the points corresponding to 
the 2-skeleton of $\Delta(t)$.
Smooth deformation
$\nu_\epsilon(t):L(t)\to M_{\Delta(t)}$ means
a smooth immersion
\[
\tilde\nu_\epsilon:L\times[0,1]
\subset M_{\tilde\Delta}
\]
such that its restriction to $M_{\Delta(t)}$
gives $\nu_\epsilon(t)$.

\begin{proof}[Proof of Corollary \ref{m-coro}]
By Proposition \ref{prop-vprod},
the mixed $h$-product coincides up to sign
with the
determinant of the linear map \eqref{eqeve}.
A small deformation of the 1-skeleton of
$\Delta(t)$
results in
a small deformation $\ll(t)$ of $\ll$.
In its turn this deformation
determines a point $q(t)\in\R^l$, $l=\kappa$,
close to the origin and such that
$h'\in\Def(h)$ passes through $\ll(t)$
iff $\ev_{\ee}(h')=q(t)$.
Since the determinant of \eqref{eqeve} is not
zero, we may find a deformation $C(t)\subset\Delta(t)$
in the class of rational even primitive tropical 
curves and lift them in the family as in
the proof of Theorem \ref{main}.
\end{proof}

\begin{exa}
The conclusion of Corollary \ref{m-coro}
may be false in the case
if $(\ez_1,\dots,\ez_\kappa)_h=0$.
For example, let 
\[
\Delta(t)=R(t)\times\R\subset\R^3,
\]
where $R(t)$ is the rectangle with vertices
$(-1,-1)$, $(-1,1)$, $(1+t,1)$, $(1+t,-1)$, $t\ge 0$.
We have 
\[
M_{R(t)}=\cp^1\times\cp^1\times\C,
\]
where the symplectic area of the first $\cp^1$
is $2\pi(2+t)$, while that of the second $\cp^1$
is $4\pi$.
The square $R(0)$ contains the interval
$C=[(-1,1),(1,-1)]$ which is an even primitive
$\Delta$-tropical curve that is Lagrangian-realizable
by an embedded sphere of homology class
$(-1,1)\in H_2(\cp^1\times\cp^1\times\C)=\Z^2$.
However the symplectic area of this class for $t>0$
is non-zero and thus it cannot be realized as 
a Lagrangian.
\begin{figure}[h]
\includegraphics[width=75mm]{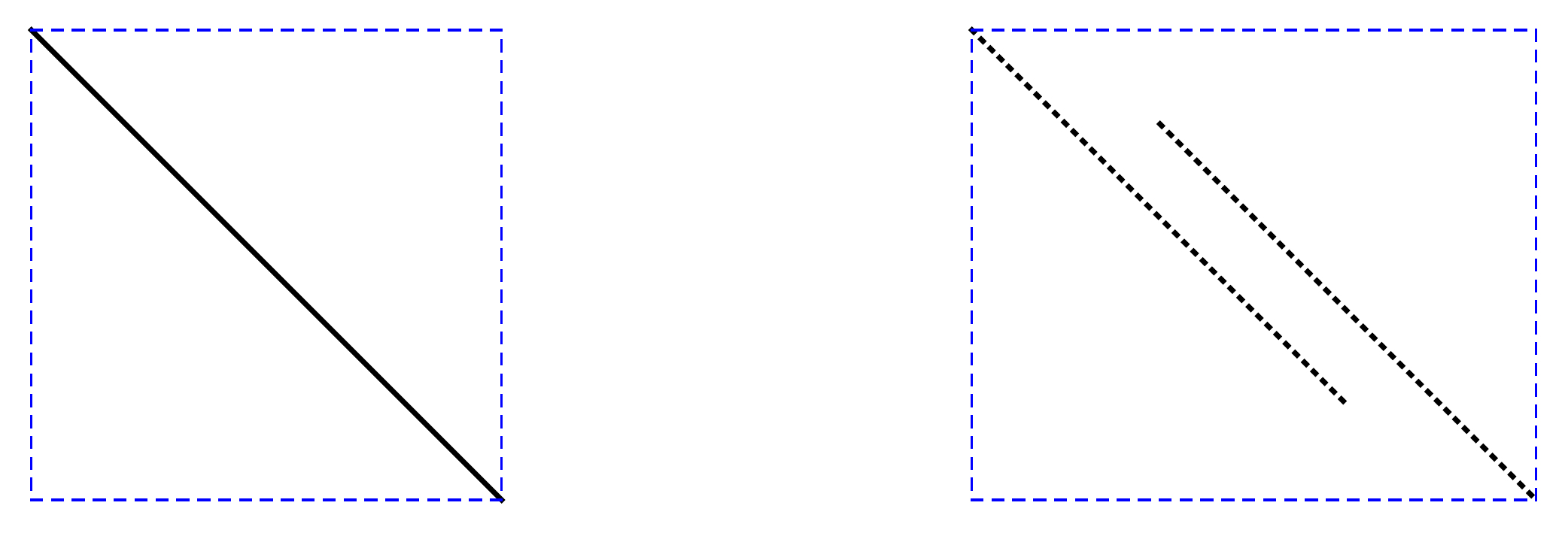}
\caption{A deformation of $\Delta$ resulting
in a disappearing tropical Lagrangian.\label{fig9}}
\end{figure}

Note that in this example $\kappa=2$,
$\ez_1=\ez_2=(1,0,0)$, so
\[
(\ez_1,\ez_2)_h=(\ez_1,\ez_2,(dh)e)=0,
\]
where $dh(e)=(1,-1,0)$ is the primitive integer
vector parallel to the only edge of $C$.
\ignore{
let $Y\subset \R^3$ be
the tropical curve obtained as the union
of three rays emanating from the origin
in the directions $(-1,-1)$, $(-1,2)$ and $(2,-1)$.
Let $\ll=\{\el_1,\el_2,\el_3\}$ be three
vertical lines passing through the points
$(-1,-1,0)$, $(-1,2,0)$ and $(2,-1,0)$.
Let $\Delta$ be the convex hull of
$\el_1\cup\el_2\cup\el_3$.
Then $Y\cap\Delta$ is an even primitive
curve in the prism $\Delta$ as all the boundary
points are bissectrice points.

However, a deformation of $\ll$ keeping $\el_2$
and $\el_3$ unchanged, and deforming $\el_1$
to a vertical line passing through $(-1,-1+t,0)$
produces a prism $\Delta(t)$  such that 
there is no Lagra
}
\end{exa}

\begin{proof}[Proof of Theorem \ref{m-thm}]
Choose a bissectrice point $p\in C\cap\dd\Delta$
as the root of the tree $C$,
and orient all edges of $C$ towards $p$.
We adopt the notation $T_e$ from the proof
of Lemma \ref{Lfeuille}. We have
$T_e=\{\iota(e)\}\times(y+\Theta_e)$.
\ignore{
for a point $\iota(e)$ inside $e$,
the 2-subtorus $\Theta_{u_e}\subset\Theta=(\R/2\pi\Z)^3$
conormal a primitive vector $u_e\in\R^3$
parallel to the edge $e\subset C$,
and a point $y\in\Theta$.
For a vector $u\in\Z^3$ not parallel to $u_e$
the intersection 
\[
\{\iota(e)\}\times\Theta_u\cap T_e
\]
is a 1-dimensional subtorus of $T_e$.
If $u\in\Z^3$ is primitive then
we denote with $S_u\subset H_1(T_e)\approx\Z$
the subgroup generated by the inclusion of
this 1-subtorus.
If $u=nu'$ for a primitive vector $u'\in\Z^3$
not parallel to $u_e$
we denote $S_u=nS_{u'}$.
}

Let
$L^+(e)\subset L$ be component
of $L\sm T_e$ that is disjoint from
$(\mdn)^{-1}(p)$.
We have $\dd L^+(e)=T_e\approx (S^1)^2$
while
\[
H_1(T_e;\R)\subset H_1( (\R/2\pi\Z)^3;\R)=\R^3
\]
is a 2-plane conormal to the edge $e$.
%
%
%
%
%
Denote $\mv(e)=\prod\limits_v m(v)$,
where the product is taken over all vertices
of $C$ that correspond to $L^+(e)$,
i.e. all the vertices that precede $e$ 
in the order corresponding to the orientation
given by the root vertex $p$.
\begin{lem}\label{lrplus}
If the rotational momentum $\rho(e)$
of an edge $e\in E_C$ is not zero then the torsion of 
$H_1(L^+(e))$ is a finite group
of order $n(e)/\mv(e)$, where $n(e)$ is the
GCD of the coordinates of $\rho(e)$,
i.e. $\rho(e)=n(e)\rho'(e)$ for a primitive vector
$\rho'(e)\in\Z^3$. 
 
Furthermore, the bivector in $\Lambda^2(\R^3)$
corresponding to $\rho(e)$ (see \eqref{bivp})
is conormal to 
the kernel of $H_1(T_e;\R)\to H_1(L^+(e);\R)$.
\end{lem}
Recall that the vector product of two vectors
$u_1,u_2\in \R^3$ may be defined  
(using the scalar product $(,)$ in $\R^3$) 
through the identity
\begin{equation}\label{bivp}
u_1\wedge u_2\wedge u = (u_1\times u_2,u)\vol
\end{equation}
that should hold for any vector $u\in\R^3$ for
the volume 3-vector $\vol\in\Lambda^3(\R^3)$
defined by the metric $(,)$.
The vector product can be thought of just as
a vector encoding
of this bivector through \eqref{bivp}
with the help of the scalar product in $\R^3$.
Since by its definition 
the rotational momentum $\rho(e)$ is obtained
as the vector product of a certain (inductively
defined) vector and $u_e$, the vectors $\rho(e)$
and $u_e$ are not parallel unless $\rho(e)=0$.
\begin{proof}[Proof of Lemma \ref{lrplus}]
Suppose that $e\subset C$ is a leaf not adjacent
to the root vertex $p$. Then $e$ is adjacent
to another bissectrice point $p_e\in C\cap\dd\Delta$.
Let $\ez_e\in\Z^3$ be a primitive tangent vector
to the edge of $\dd\Delta$ containing the point $p_e$. 
By the construction of $M_\Delta$,
$L^+(e)\approx S^1\times D^2$ and
the kernel of $H_1(T_e;\R)\to H_1(L^+(e);\R)$
is cut by the conormal torus of $\ez_e$
(i.e. the torus conormal to the bivector encoded
by $\rho(e)$ through \eqref{bivp}).
In this case there is no torsion
in $H_1(L^+(e);\Z)\approx\Z$.
Also $\rho(e)$
is primitive as $p_e$ is a bissectrice point.

Let $e$ be an oriented edge
whose source is a 3-valent vertex $v\in V_C$.
By the induction hypothesis we know that the lemma
already holds for the two edges
$\iota_1,\iota_2$ that are incoming with
respect to $v$.
Suppose at first that $m(v)=1$. Then $\nu_\epsilon$
induces an isomorphism between $H_1(Q_v)$
and $H_1(\C^3)=H_1(\Theta)=\Z^3$ for 
the component $Q_v\approx P\times S^1$
of $L\sm T$ corresponding
to $v$ (i.e. bounded by the tori $T_{\iota_1}$,
$T_{\iota_2}$ and $T_e$).
The kernel $K_e$ of $H_1(T_e)\to H_1(L^+(e))$
is cut in $H_1(T_e)\subset \Z^3$ 
by $K_{\iota_1}+K_{\iota_2}$ (which is two-dimensional
by the induction hypothesis unless $\rho(e)=0$).


Let $\Lambda_\iota\subset\Z^3$ be
the lattice formed by
the integer vectors 
parallel to the elements of $K_{\iota_1}+K_{\iota_2}$.
Consider the inclusion 
homomorphism \[
H_1(L^+(\iota_1)\cup L^+(\iota_2))\to
H_1(L^+(e))
\]
from the exact sequence of the pair
$(L^+(e),L^+(\iota_1)\cup L^+(\iota_2))$.
By the excision property,
the homology groups to the left and to the right
of this homomorphism in the long exact sequence are torsion-free.
Thus the torsion of $H_1(L^+(e))$ decomposes
into the sum 
of the torsions of $H_1(L^+(\iota_j))$, $j=1,2$,
and the quotient
$\Lambda_\iota/(K_{\iota_1}+K_{\iota_2})$.
Thus $\#(H_1(L^+(e)))$ equals to
$n(\iota_1)n(\iota_2)$
times the GCD of the coordinates of
$\rho'(\iota_1)\times\rho'(\iota_2)$.
The kernel $K_e\otimes\R$ is the intersection
of $K_{\iota_1}\otimes\R+K_{\iota_2}\otimes\R$
with $H_1(T_e;\R)\subset H_1(\Theta;\R)=\R^3$
and thus is conormal to the bivector
corresponding to 
$\rho(\iota_1)\times\rho(\iota_2)$.

If $m(v)>1$ then the homomorphism
$H_1(Q_v)\to H_1(\C^3)=\Z^3$ 
induced by 
$\nu_\epsilon$ is injective, but not surjective.
Its image is a sublattice of index $m(v)$.
Thus to get $\#(H_1(L^+(e)))$ we have to divide
the product of  $n(\iota_1)n(\iota_2)$ and
the GCD of the coordinates of
$\rho'(\iota_1)\times\rho'(\iota_2)$ by $m(v)$.
The kernel
of $H_1(T_e;\R)\to H_1(L^+(e);\R)$
is defined over $\R$, thus its computation
is the same as in
the $m(v)=1$ case.
%
\ignore{
we can lift $Q_v$ to 
a multiplicative linear map
$(\C^\times)^3\to (\C^\times)^3$ 
whose (multiplicative) matrix is made of
the coordinates
of the primitive vectors tangent to $\iota_1$
and $\iota_2$ and $\frac1{m(v)}\iota_1\times\iota_2$.
It is a covering of degree $m(v)$ 
and the (cf. \cite{Mi-05}).

While the topology upstairs is the 

By the induction hypothesis, $\rho$ is conormal to
the linear span of 
}
\end{proof}
To finish the proof of Theorem \ref{m-thm} we
apply Lemma \ref{lrplus} to the leaf $e_p$ of $C$
adjacent to the bissectrice point $p$.
The manifold $L$ is obtained
by gluing of $L^+(e_p)$ and $S^1\times D^2$
along their common boundary $T_e\approx (S^1)^2$.
The kernels of the maps from $H_1(T_e)$ into
these 3-folds are isomorphic to $\Z$
and described by the rotational momenta $\rho(e)$
and $\rho_p$.
The cardinality $\#(H_1(L))$ 
equals to the product of the cardinality
of the torsion subgroup of
$H_1(L^+(e_p))$ and
the GCD of the coordinates of $\rho'(e)\times\rho_p$
(recall that $\rho_p$ is primitive as $p$ is
a bissectrice point), i.e. to
$(\ez_1,\dots,\ez_\kappa)_h$ in the case when the
latter quantity is non-zero.
If $(\ez_1,\dots,\ez_\kappa)_h=0$ then 
$H_1(L;\R)\neq 0$.
\end{proof}

\begin{rmk}
Note that the topology of $L$ is determined
by the tropical curve $C\subset L$ and the lines
containing the edges of $\Delta$ incident
to $C$. In other words, it is determined by
the tropical curve $h:\Gamma\to\R^3$ and
the lines $\ll$ forming a boundary configuration
for $h$. Furthermore, we may construct the 3-manifold $L$
even in the case when no suitable Delzant domain
for $(h,\ll)$ exists. 

Upgrading the tropical multiplicity 
$m(h,\ll)$ to a torsion group $H_1(L)$ 
of order $m(h,\ll)$ may be considered as
a certain refinement of $m(h,\ll)$.
It might be interesting to attempt to extract
finer invariants in tropical enumerative
problems using $H_1(L)$, and perhaps even
more interesting using other invariants
of $\pi_1(L)$
or
the graph 3-fold $L$ itself.
\end{rmk}

\begin{defn}
Let $h:\Gamma\to\R^3$ be a primitive tropical immersion
of a graph $\Gamma$ with $\kappa$ ends,
and $\ll=(\el_1,\dots,\el_\kappa)$
be a generic boundary configuration of lines for $h$
parallel
to $\zz=(\ez_1,\dots,\ez_\kappa)$.
We say that a Delzant polyhedral domain
$\Delta\subset\R^3$
is {\em $(h,\ll)$-suitable}
if for each $j=1,\dots,\kappa$
there exists an edge $e_j$ of $\Delta$ parallel
to $\ez_j$ such
that $e_j\supset h(\Gamma)\cap\el_j$ and all points
of $\dd\Delta\cap h(\Gamma)$ are bissectrice
points.
\end{defn}

Determining existence of a $(h,\ll)$-suitable
domain seems to be a non-trivial problem.
Still it is easy to establish some necessary conditions.
Let $\DD=(d_1,\dots,d_\kappa)$ be the toric degree
of $h$ and $p_j\in\el_j$ be the points of intersection
of the line $\el_j$ and the leaf of $h(\Gamma)$ 
parallel to $d_j\in\Z^3$.
Let $P$ be the convex hull of
$\bigcup_{j=1}^\kappa\{p_j\}$.
\begin{prop}\label{ncondD}
If $\Delta\subset\R^3$ is an $(h,\ll)$-suitable 
polyhedral domain for a generic configuration
$\ll$ then $d_j\times\ez_j\in\Z^3$
is primitive and 
$p_j$ is a vertex of the polyhedron $P$
for every $j=1,\dots,\kappa$. 
\end{prop}
\begin{proof}
Since an $(h,\ll)$-suitable 
polyhedral domain $\Delta$ is convex and $p_j$
is contained in the 1-skeleton of $\dd\Delta$,
the point $p_j$ must be contained in the 1-skeleton
of $P$. If $p_j$ is not a vertex then it is contained
in the interval between two other bissectrice
points all sharing the same line from $\ll$.
This contradicts to the assumption that $\ll$
is generic.

Since $p_j$ is a bissectrice point for an 
edge of the Delzant domain $\Delta$
parallel to $\ez_j$, the vector $d_j$ can be 
presented as the sum of two vectors $a_j$ and $b_j$
(parallel to the adjacent faces of $\Delta$)
such that $a_j,b_j,\ez_j$ is a basis of $\Z^3$.
Thus $d_j\times\ez_j$ is primitive.
\end{proof}

The first necessary condition given by Proposition
\ref{ncondD} is also sufficient to present the
compact 3-manifold $L$ corresponding to $h$
as a smooth Lagrangian manifold in a non-compact
symplectic toric manifold corresponding to,
perhaps, a non-convex domain $\Delta$
(but still convex and polyhedral 
near its boundary faces),
cf. Remark \ref{nonconvD}. 

Let $\ll$ be a generic boundary configuration for
a primitive tropical curve $h:\Gamma\to\R^3$.
Denote with $C\subset\R^3$ the closure of the
the only bounded component of
$h(\Gamma)\setminus \bigcup\limits_{j=1}^\kappa(p_j)$.
Note that for any $(h,\ll)$-suitable Delzant
domain $\Delta$ we have $C=h(\Gamma)\cap\Delta$.

\begin{prop}\label{nclift}
If $\ll$ is a generic boundary configuration
for a primitive tropical immersion $h:\Gamma\to\R^3$,
and $d_j\times\ez_j\in\Z^3$ is primitive for any
$j=1,\dots,\kappa$
then there exists a set $\Delta\subset\R^3$
and its subset $\dd\Delta\subset\R^3$ such that
the following conditions hold.
\begin{enumerate}
\item $\Delta\supset C$, $p_j\in\dd\Delta$, $j=1,\dots,
\kappa$.
\item 
$\Delta$ is locally a Delzant polyhedral domain
near any point of $\dd\Delta$ and $p_j$ is a bissectrice
point so that $C$ is an even primitive tropical 
curve in $\Delta$.
\item The set $\Delta\setminus\dd\Delta$ is open.
\item The curve $C\subset\Delta$ is Lagrangian-realizable
in $M_{\Delta}$ in the sense of Definition \ref{dLr}.
\end{enumerate}
\end{prop}
Thus we may speak of a Lagrangian graph manifold
$L$ for $(h,\ll)$ even in the absence of an
$(h,\ll)$-suitable Delzant domain in $\R^3$.
\begin{proof}
Since $d_j\times\ez_j$ is primitive,
we may find $a_j,b_j\in\Z^3$ such that $a_j+b_j=-d_j$,
and $(a_j,b_j,\ez_j)=1$, i.e. such that
$a_j,b_j,\ez_j$ is a basis of $\Z^3$.
Define $\Delta_{a_j}$ (resp. $\Delta_{b_j}$)
to be 
the half-space passing through $p_j$
whose boundary is parallel
to $\ez_j$ and $a_j$ (resp. $b_j$) and containing
the edge of $C$ adjacent to $p_j$.
The intersection $\Delta_j=\Delta_{a_j}\cap
\Delta_{b_j}$ is a Delzant domain which
has $\el_j$ as its apex edge.
We set
\[
\Delta=U_C\cup\bigcup\limits_{j=1}^\kappa
(U_{p_j}\cap\Delta_j),
\]
where $U_{p_j}$ is a ball around $p_j$ of a small
radius, and $U_C$ is a very small open
neighborhood of $C$.
The proof of Theorem \ref{main} produces
Lagrangian realizability of $C\subset\Delta$
in the (non-compact) symplectic manifold $M_\Delta$
in the same way as in the case when $\Delta$ is a
global Delzant domain.
\end{proof}

The following two examples show that
the topology of the Lagrangian 3-fold $L$
is not determined by the corresponding
tropical multiplicity even in the case
when the multiplicity is 1.
In the first example
we may find an $(h,\ll)$-suitable Delzant domain $\Delta$.
In the second example we do not claim
existence of an $(h,\ll)$-suitable Delzant domain,
though we can still easily lift $L_8$
to 
a singular
toric 3-fold $M_\Delta$ (for non-Delzant compact
polyhedral
domain $\Delta$) so that $L_8$ intersect the 
singular locus of $M_\Delta$. Also Proposition
\ref{nclift} provides a Lagrangian embedding
of $L_8$ to a non-compact toric symplectic manifold.

\begin{exa}[Lens spaces]\label{exa-ls}
Suppose that 
\[
C=[(0,0,0),(0,0,1)]\subset\R^3
\]
is the interval 
,
and $\Delta_{p,q}\supset C$ is a Delzant polyhedral domain
such that it possesses edges contained
in the line $\el_1$ through $(0,0,0)$ parallel
to the vector $\ez_1=(1,0,0)$ and 
in the line $\el_2$ through $(0,0,1)$ parallel
to the vector $\ez_2=(-q,p,0)$ for an integer $p>0$
and an integer $q$ relatively prime with $p$,
and such that $(0,0,0)$ and $(0,0,1)$
are bissectrice points of $C$.
Then Theorem \ref{main} gives an embedded Lagrangian
$L_{p,q}\subset M_{\Delta_{p,q}}$ diffeomorphic to the
lens space $L(p,q)$. Note that the tropical multiplicity
is $p$ while $L(p,q)$ and $L(p,q')$ may be
non-diffeomorphic.
\begin{figure}[h]
\includegraphics[width=45mm]{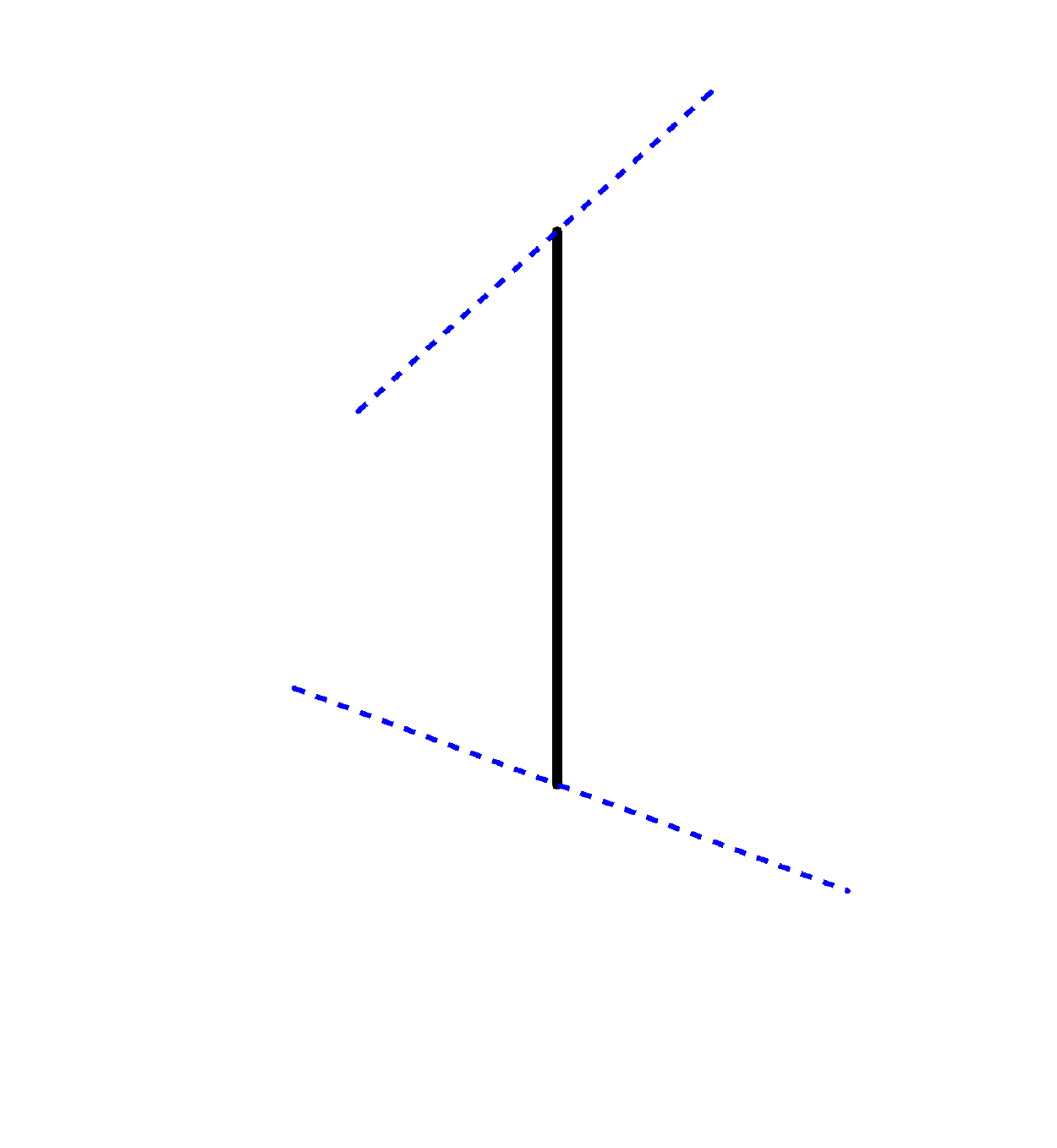}
\caption{A Lagrangian lens space. \label{fig10}}
\end{figure}

It is easy to see that such $\Delta_{p,q}$ exists
for arbitrary $p,q$.
Since $(0,0,1)\times \ez_j\in\Z^3$, $j=1,2$,
is primitive, we may present $\el_j$ as
the apex edge of the intersection of
two half-space such that
$C$ has the boundary momentum 1 with respect
to each of them. The intersection of the four
half-space is a tetrahedron $\Delta'_{p,q}$ which
has rational slope, but is, perhaps, non-Delzant
at some of its faces disjoint from $C$.
To get $\Delta_{p,q}$ we truncate
the tetrahedron $\Delta'_{p,q}$ at
such faces. Note that a truncation
$\Delta_{p,q}\subset\Delta'_{p,q}$ may be associated
to a toric
resolution $M_{\Delta_{p,q}}\to M_{\Delta'_{p,q}}$ of the
toric orbifold $M_{\Delta'{p,q}}$.
\end{exa}

\begin{exa}[The Poincar\'e sphere]
\label{exa8}
Let
$e_1=[(-1,0,0),(0,0,0)]$, $e_2=[(0,-1,0),(0,0,0)]$,
$e_3=[(0,0,0),(1,1,0)]$,
\[ 
C=e_1\cup e_2\cup e_3,
\]
and $\el_j$, $j=1,2,3$,
be the lines passing through $p_1=(-1,0,0)$,
$p_2=(0,-1,0)$, $p_3=(1,1,0)$, and parallel
to $\ez_1=(0,1,2)$, $\ez_2=(1,0,3)$,
$\ez_3=(0,1,5)$. 
The corresponding 3-manifold $L_8$
is obtained by gluing three solid tori to
the product $P\times S^1$
of the pair-of-pants $P$ and the circle $S^1$
according to the rotational momenta of the leaves of $C$
which are $\rho_1=(0,2,-1)$, $\rho_2=(-3,0,1)$ and
$\rho_3=(5,-5,1)$.

It is easy to see
that $L_8$ is the Poincar\'e homology sphere,
i.e. the Seifert-fibered homology sphere with 
three multiple fibers of Seifert invariant
$(2,1)=(2,-1)$, $(3,1)$ and $(5,1)$,
cf. \cite{KiSch}.
\begin{figure}[h]
\includegraphics[width=80mm]{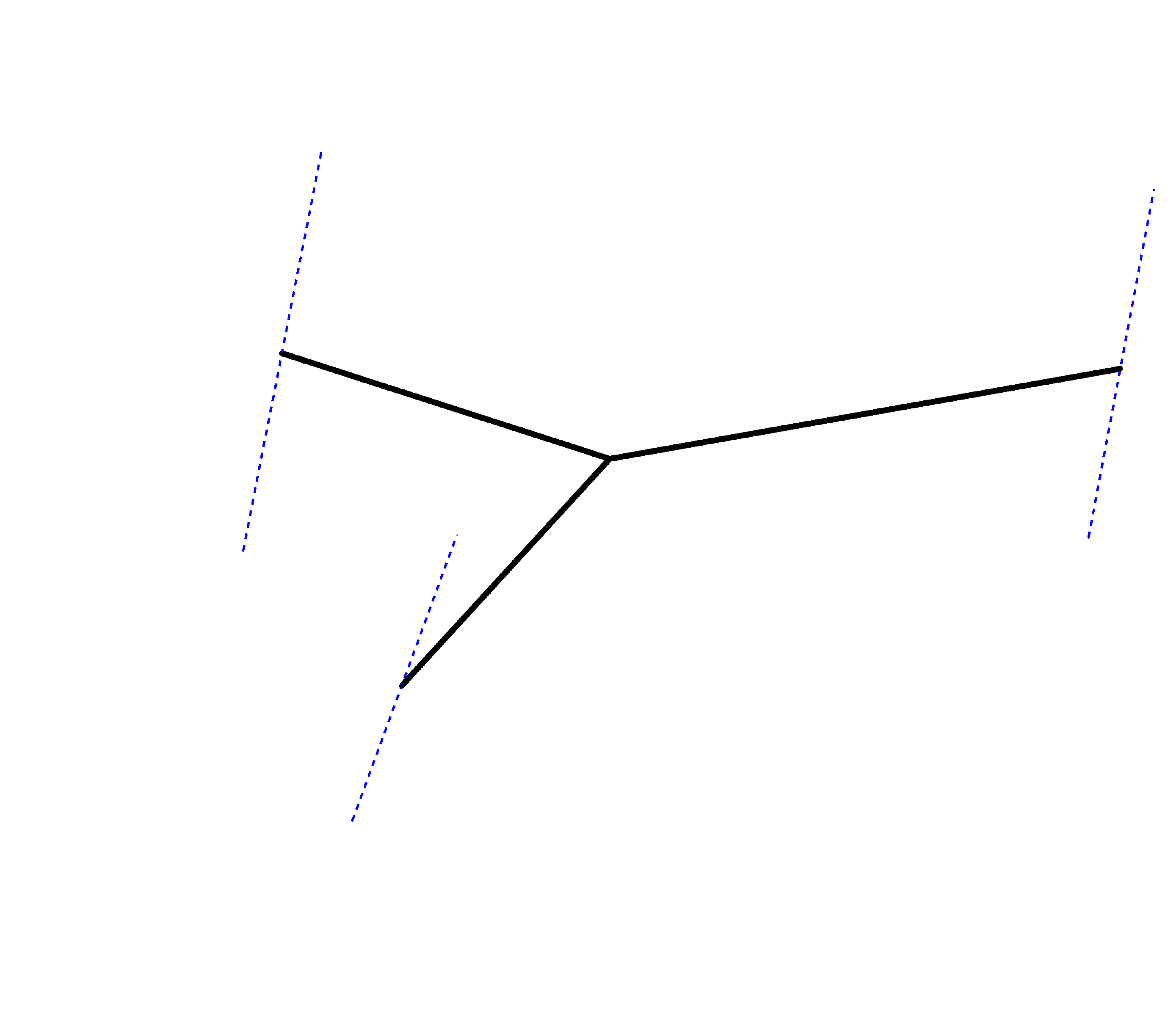}
\caption{The Poincar\'e sphere from a tropical curve.
\label{fig11}}
\end{figure}

We may find a compact 
polyhedral domain $\Delta_8\supset C$
with edges parallel to $\ez_j$, containing $e_j$.
A straightforward modification of
Theorem \ref{main} produces
a Lagrangian mapping of $L_8$
to the orbifold $M_{\Delta_8}$.

To find
$\Delta_8$
we
define the convex domains $\Delta_j$, $j=1,2,3$,
to be the convex hull
of $\el_j$ and the rays emanating from $p_j$
in the direction of $(1,1,0)$ and $(0,-1,0)$
for $j=1$; $(-1,0,0)$ and $(1,1,0)$ for $j=2$;
$(0,-1,0)$ and $(-1,0,0)$ for $j=3$, see Figure
\ref{fig11}.
We set
\[
\Delta_8=\{|x|\le \epsilon\}\cup \bigcup_{j=1}^3\Delta_j
\]
for a small $\epsilon>0$. It is a (non-Delzant) polyhedron
with eight facets.

Perturbing $\Delta_8$ by introducing more faces to resolve
singularities of $M_{\Delta_8}$ we may obtain a Delzant polyhedron
$\tilde\Delta_8$ also containing $C$ in a way that the endpoints of $C$
sit on the edges of $\tilde\Delta_8$ parallel to $\el_j$.
This gives us a Lagrangian mapping of $L_8$ to the
smooth symplectic toric manifold $M_{\tilde\Delta_8}$.
Note, however, that this mapping is not an embedding
unless the edges $e_j$ can be made
bissectrices of the corresponding
angles. 
\end{exa}

\begin{que}
Does $L_8$ admit a Lagrangian embedding to a
smooth symplectic toric manifold?
\end{que}

We have $m(h,\ll)=1$ for the associated tropical problem
both in Example \ref{exa8} and in Example \ref{exa-ls}
for $p=1$.
Nevertheless, topology of the corresponding
Lagrangian manifolds is different: $L_{1,0}$ is
the standard sphere $S^3$ while
$L_8$ is the (non simply connected) Poincar\'e sphere.

\begin{exa}
Let $\Delta$ be the convex hull of $(0,0,0)$,
$(1,0,0)$, $(0,1,0)$ and $(0,0,1)$.
The baricenter of $\Delta$ is $p=(\frac14,\frac14,\frac14)$.
Let \[
C_{12}=[(\frac14,0,0),p]\cup [(\frac12,\frac12,0),p]\cup 
[(0,\frac14,\frac34),p].\]
It is easy to see that $C_{12}$
is a primitive tropical curve in $\Delta$ while 
the multiplicity of its only 3-valent vertex $p$ is 1.
Thus there exists a Lagrangian embedding 
\[
L_{12}\subset M_\Delta=\cp^3
\]
for a rational homology sphere $L_{12}$ with
\[
\#(H_1(L_{12}))=|(\ez_1,\ez_2,\ez_3)_h|=4,
\]
where $\ez_1=(1,0,0)$, $\ez_2=(1,-1,0)$,
$\ez_3=(0,1,-1)$ are primitive integer vectors in
the directions parallel to the edges of $\Delta$
containing the endpoints of $C_{12}$.

It is interesting to compare the resulting $L_{12}$ against
Chiang's example of a nonstandard Lagrangian
$L_C\subset\cp^3$, see \cite{Chiang}.
It is easy to see $\pi_1(L_{12})=\pi_1(L_C)$
(in particular, that 
the fundamental group  
of $L_{12}$ has order 12).
\end{exa}

\begin{que}
Is $L_{12}$ Hamiltonian isotopic to Chiang's Lagrangian submanifold
$L_C\subset\cp^3$?
\end{que}

\begin{que}
Are there other rational homology 3-spheres (except for $\rp^3$
and $L_{C}$) that admit Lagrangian embeddings to $\cp^3$?
\end{que}

Note that according to Seidel's theorem \cite{Seidel-L}
we have $\#(H_1(L))\equiv 0\pmod2$
for any Lagrangian rational
homology sphere $L\subset\cp^3$.